\documentclass[reqno,11pt]{amsart}

\usepackage{amsthm}
\usepackage{amssymb}
\usepackage{graphicx}			
\usepackage{appendix}
\usepackage{mathrsfs}

\usepackage{fontenc}
\usepackage{inputenc}

\usepackage{hyperref}

\usepackage{tikz-cd}

\usepackage{enumitem}

\usepackage{mathtools}

\newtheorem{thm}{Theorem}[section]
\newtheorem{lem}[thm]{Lemma}
\newtheorem{cor}[thm]{Corollary}
\newtheorem{prp}[thm]{Proposition}

\theoremstyle{definition}
\newtheorem{dfn}[thm]{Definition}
\newtheorem{qst}[thm]{Question}

\theoremstyle{remark}
\newtheorem{rem}[thm]{Remark}

\newcommand{\R}{\mathbb{R}}
\newcommand{\N}{\mathbb{N}}
\newcommand{\Q}{\mathbb{Q}}
\newcommand{\Z}{\mathbb{Z}}

\newcommand{\Orb}{\mathcal{O}}
\newcommand{\RP}{\mathbb{RP}}
\newcommand{\Sph}{\mathbb{S}}

\newcommand{\PT}{\mathbb{P}T^*}
\newcommand{\ST}{\mathbb{S}T^*}

\usepackage[nodayofweek, level]{datetime}

\title[Besse families of hypersurfaces]{Rigidity and deformations of Besse \\ families of minimal hypersurfaces}

\author[Christian Lange]{Christian Lange}
\address[Christian Lange ]{Ludwig-Maximilians-Universit\"at M\"unchen, Mathematisches Institut\newline\indent Theresienstr. 39, 80333 München, Germany}
\email{lange@math.lmu.de}


\date{2nd December 2026}

\begin{document}

\maketitle

\begin{abstract}
We generalize Zoll families of hypersurfaces to Besse families of hypersurfaces, which are defined differential-topologically and may contain non-embedded and exceptional hypersurfaces. A contact-geometric discussion presents them as double Legendrian fibrations and shows that regular hypersurfaces are either spheres or real projective spaces.
As an analogue of Wadsley's theorem for Riemannian metrics all of whose geodesics are closed, we show that a family of closed immersed minimal hypersurfaces whose Legendrian lifts foliate the space of contact elements is a Besse family. In analogy with the Berger conjecture on simply connected manifolds all of whose geodesics are closed, we show that every hypersurface of a Besse family on a closed manifold is regular. On even-dimensional manifolds we obtain a stronger conclusion, generalizing a theorem of Gromoll and Grove about $2$-spheres all of whose geodesics are closed: every hypersurface of a Besse family is regular and embedded. Moreover, building on work by Ambrozio and Guajardo, we provide metrics with trivial isometry group carrying Zoll families of non-embedded minimal hypersurfaces in every odd dimension at least three, and Besse families of minimal hypersurfaces with exceptional hypersurfaces on orbifolds in every dimension at least two. We also give a short proof that manifolds admitting Zoll families of embedded hypersurfaces are homeomorphic to spheres or real projective spaces, recovering a recent result of Martins up to homeomorphism. In particular, this applies to all even-dimensional manifolds with a Besse family.
\end{abstract}

\section{Introduction}

Following \cite{MS22}, we call a Riemannian metric \emph{Besse} if all its geodesics are closed, \emph{Zoll} if in addition all prime geodesics have the same length, and \emph{simple Zoll} if further all of them are embedded. See \cite{Zo03} for the first nontrivial examples, \cite{Gu76,LM02,AG26} for more sophisticated constructions, and \cite{Be78} for a foundational overview. Ambrozio, Marques and Neves \cite{AMN25} generalized the notion of a simple Zoll Riemannian surface to the notion of a Zoll family of embedded minimal hypersurfaces and constructed nontrivial deformations starting from the model example of the round sphere $\mathbb{S}^{n+1}$, $n\geq2$, with its equatorial hyperspheres. Ambrozio and Guajardo afterwards obtained similar deformations of the real projective space $\mathbb{RP}^{n+1}$, $n\geq2$ \cite[Theorem~1.5]{AG26}. On $\mathbb{S}^3$~\cite{AMN24} and on $\mathbb{RP}^3$~\cite{Ma26}, Riemannian metrics admitting such Zoll families are characterized by the equality of certain area widths.

Martins further abstracted this notion in a differential-topological way, introducing \emph{unoriented transitive families of hypersurfaces}, which are a priori not required to be embedded, and proved that, for $n\geq 2$, an $(n+1)$-manifold carrying such a family of closed embedded hypersurfaces is diffeomorphic to $\mathbb{S}^{n+1}$ or $\mathbb{RP}^{n+1}$, that the hypersurfaces are diffeomorphic to $\mathbb{S}^{n}$ resp.\ $\mathbb{RP}^{n}$, and that the parameter space is homeomorphic resp.\ diffeomorphic to  $\mathbb{RP}^{n+1}$ \cite{Ma26}. We independently developed shorter arguments which, together with the proof of \cite[Proposition~3.5]{Ma26}, yield the classification up to homeomorphism for all $n\geq 1$; see Section~\ref{sec:topo_class}. Both Martins' approach and ours rest on the idea of isotoping a two-sided hypersurface of the family onto itself with reversed coorientation. This is achieved by rotating it, through other members of the family, about codimension two submanifolds; cf.\ \cite{BB77} and \cite{Be78} for an analogous argument for geodesics. The one-sided case reduces to this by passing to a double cover. Our proof is given in Section~\ref{sec:topo_class}, where such \emph{isotopically reversible hypersurfaces} are studied in a greater generality than needed here, which may be of independent interest.

The conclusion of the topological classification fails without the embeddedness assumption. Indeed, the family of equatorial hyperspheres on a  round odd-dimensional sphere $\mathbb{S}^{2m+1}$, $m\geq1$, descends to Zoll families of non-embedded minimal hypersurfaces on certain odd-dimensional spherical space forms, e.g. lens spaces. Moreover, these examples can be deformed. The equatorial hyperspheres on $\mathbb{S}^{2m+1}$ are also minimal with respect to the Berger metrics, which are obtained from the round metric by rescaling along the fibers of the Hopf fibration~\cite{To12,Am26}. Since the Berger metrics are invariant under the scalar action of $\mathrm{U}(1)$, the lens spaces $\mathbb{S}^{2m+1}/\mathbb{Z}_k$ for the cyclic subgroups $\mathbb{Z}_k<\mathrm{U}(1)$ of order $k\geq3$ carry Zoll families of non-embedded minimal hypersurfaces with respect to the non-round Berger metrics, which have non-constant curvature. In addition to that, building on the machinery of Ambrozio and Guajardo \cite{AG26}, who constructed Zoll metrics on the $2$-sphere with prescribed isometry group \cite[Theorem 1.4]{AG26}, we obtain the following statement valid in all dimensions.

\begin{thm}\label{thm:zoll_finite_isometry} For $n\geq 1$ let $G < \mathrm{O}(n+2)$ be a finite subgroup with $-\mathrm{id} \notin G$. Then there are embedded spheres in $\R^{n+2} \setminus \{0\}$ whose induced metrics admit a Zoll family of embedded minimal hyperspheres, whose isometry groups are $G$ and such that $G$ preserves the Zoll family.
\end{thm}

The assumption  $-\mathrm{id} \notin G$ is necessary for $n=1$ as every Riemannian metric with all geodesics closed on the real projective plane has constant curvature \cite{Pr09,Gr63}.

The quotients of the Berger spheres are still highly symmetric. In contrast, combining Theorem~\ref{thm:zoll_finite_isometry}, applied to free actions of finite groups not containing $-\mathrm{id}$, with Proposition~\ref{Berger_odd_orientable}, we obtain examples without any symmetries.

\begin{thm}\label{thm:non_embed_odd}
In every odd dimension at least three, there exist Riemannian manifolds with trivial isometry group carrying a Zoll family of non-embedded minimal hypersurfaces.
\end{thm}

The construction above applies equally if the finite group does not act freely. The quotient is then an orbifold, and the Legendrian lifts of the images of the equatorial hyperspheres only define an orbifold fiber bundle. Similarly, by Wadsley's theorem \cite{Wa75} the geodesic flow of a Besse Riemannian metric on the unit sphere bundle factors through an almost free circle action and hence defines an orbifold circle bundle over the space of closed geodesics. We define an \emph{unoriented Besse transitive family of hypersurfaces} on a smooth $(n+1)$-manifold $M$ with $n\geq1$ as a collection $\mathscr B$ of immersed, closed, connected hypersurfaces $f\colon\Sigma\to M$, considered up to reparametrization, whose Legendrian lifts $\hat f\colon\Sigma\to\PT M$, $p\mapsto T_pf(T_p\Sigma)$, to the projectivized cotangent bundle $\mu\colon\PT M\to M$ are embeddings with images precisely the fibers of a smooth orbifold fiber bundle $\nu\colon\PT M\to P$. Here we identify a projectivized covector with its kernel, so that points of $\PT M$ are tangent hyperplanes of $M$; see Section~\ref{sec:contact_geometry} for a contact-geometric description. We refer to such a family as a \emph{Besse family} for short. We call a hypersurface of a Besse family \emph{regular} if it corresponds to a manifold point of the orbifold $P$ and \emph{exceptional} otherwise. The unoriented transitive families of closed hypersurfaces of \cite{Ma26} are exactly those Besse families for which $P$ is a manifold, i.e.\ for which all hypersurfaces of the family are regular. In this case we refer to the family as a \emph{Zoll family}. In contrast, in \cite{AMN25, AG26,Ma26} the term ``Zoll'' signifies that the hypersurfaces in the family are closed, embedded and minimal. Apart from that, our notation largely follows that of \cite{Ma26} to facilitate comparison.

A Besse $2$-orbifold in the sense of \cite{La20} is a $2$-orbifold with a Besse Riemannian metric. This notion and the examples arising from Theorem~\ref{thm:zoll_finite_isometry} as quotients of non-free actions can be subsumed by the notion of a Besse family of minimal hypersurfaces by also allowing the ambient manifold $M$ and the hypersurfaces in the definition to be orbifolds. Rotationally symmetric Besse $2$-orbifolds already appeared in the work of Tannery \cite{Ta92}. Non-rotationally symmetric examples can be obtained from \cite[Theorem~1.4]{AG26} and by the approach in \cite{LM02}. An argument similar to the one in the proof of Theorem~\ref{thm:non_embed_odd}, but applied to non-free group actions, yields the following consequence of Theorem~\ref{thm:zoll_finite_isometry}.

\begin{cor}\label{cor:Besse_orbifold} For every $n\geq 1$ there exist closed Riemannian $(n+1)$-orbifolds of non-constant curvature which admit a non-Zoll Besse family of minimal hypersurfaces.
\end{cor}

For the family to be non-Zoll it suffices that the orbifold has non-isolated singularities; see Remark~\ref{rem:non-Zoll}.

The question arises whether there exist non-Zoll Besse families of hypersurfaces on \emph{manifolds}. In the geodesic setting, the answer is affirmative. The round metric on an odd-dimensional sphere descends to non-Zoll Besse metrics on certain quotients, e.g.\ lens spaces. However, Berger conjectured that a simply connected Riemannian manifold all of whose geodesics are closed is Zoll, i.e.\ that all its geodesics have the same length; cf.\ \cite{Be78}. This statement was first proved for $2$-spheres by Gromoll and Grove \cite{GG81}, and is by now a theorem for all spheres \cite{RW17,Di25}, see also \cite{Di26a,Di26b}. For Besse families of hypersurfaces, we show the following purely differential-topological statement.

\begin{thm}\label{thm:Berger}
A Besse family of hypersurfaces on a closed manifold is Zoll, i.e.\ all its hypersurfaces are regular.
\end{thm}

In particular, in contrast to the geodesic setting no simple connectedness assumption is needed. On the other hand, the manifold assumption is essential, since non-Zoll Besse families exist on orbifolds in every dimension; see Corollary~\ref{cor:Besse_orbifold}. 

As for the topological classification of Zoll families of hypersurfaces, a first key observation for the proof of Theorem~\ref{thm:Berger} is that all regular hypersurfaces of a Besse family of closed hypersurfaces must be spheres or real projective spaces, and in the even-dimensional case in fact spheres; see Lemma~\ref{lem:lifted_orientations}. 
This observation is a consequence of  well-known facts in contact geometry; see Subsection~\ref{sub:coor_and_lifts}. Related contact-geometric arguments identify the cosphere bundle of $M$ with that of the parameter space $P^+$ of the cooriented hypersurfaces of the family, and translate fiberwise involutions and covering actions into actions on $P^+$. This quickly settles the orientable odd-dimensional case; see Proposition~\ref{Berger_odd_orientable}. In the orientable even-dimensional case, Gysin sequences and Poincaré duality show that $M$ and $P^+$ are rational homology spheres; see Proposition~\ref{prp:even_rat_hom_sphere}. Singular points are then excluded by Satake's orbifold Gauß--Bonnet theorem~\cite{Sa57}, combined with a Mayer--Vietoris argument in Borel cohomology, and in dimension two by the classification of Seifert fibrations; see Proposition~\ref{prp:berger_even_dim_orient}. The non-orientable case is reduced to the orientable one via the orientable double cover and the Lefschetz fixed point theorem; see Subsection~\ref{sub:Berger_non_orientable}.

For a Hamiltonian generalization of the statement that all geodesics on a Besse $2$-sphere have the same length, see \cite{FLS19}. For \emph{oriented} Besse families, defined analogously with cooriented hypersurfaces and the bundle of cooriented tangent hyperplanes, the conclusion of the theorem already fails in dimension two. Counterexamples of geometric origin are given by \emph{irreversible} Katok Finsler $2$-spheres with rational parameter, all of whose geodesics are closed \cite{Ka73}.

The theorem of Gromoll and Grove not only says that all prime geodesics of a Riemannian metric on $\Sph^2$ all of whose geodesics are closed have the same length, but also that all of them are embedded. For Zoll $2$-spheres their proof is based on the existence of a simple closed geodesic and an open-closed connectedness argument; cf. Lemma~\ref{lem:simple_charact}. The following  differential-topological statement holds on all even-dimensional manifolds and does not rely on the existence of an embedded hypersurface in the family.

\begin{thm}\label{thm:Gromoll--Grove} An unoriented transitive family of closed hypersurfaces on an even-dimensional manifold $M$ is simple Zoll, i.e.\ all its hypersurfaces are embedded. Moreover, $M$ is closed and the hypersurfaces in the family are two-sided if and only if $M$ is orientable. 
\end{thm}

The proof of Theorem~\ref{thm:Gromoll--Grove} is based on a hairy-ball-type obstruction. A double point of a hypersurface yields a section of the pullback of the sphere bundle $\nu^+\colon ST^*M\to P^+$ to the manifold of ordered double points, which fibers over $P^+$ with fibers of dimension less than $n$. The fact that $P^+$ is a rational homology sphere whose fundamental group has no proper subgroups of finite index (see Proposition~\ref{prp:even_rat_hom_sphere} and Corollary~\ref{cor:no_finite_covers}), implies that this pullback is injective in top-degree rational cohomology, which contradicts the nonvanishing of the Euler class of $\nu^+$; see Propositions~\ref{prop:I} and~\ref{prp:Gromoll_Grove_even_dimensional}. The non-orientable case reduces to the orientable one via the orientable double cover; see Proposition~\ref{prp:Gromoll_Grove_non-orien}.

Combining Theorems~\ref{thm:Berger} and~\ref{thm:Gromoll--Grove} with Proposition~\ref{prp:besse_impliese_cmpact} and Martins' classification~\cite{Ma26} (resp.\ Theorem~\ref{thm:top_class} and~\cite{Th97} for $n=1$), we obtain the following generalization of the Gromoll--Grove theorem together with a topological classification.

\begin{cor}\label{cor:Gromoll--Grove}
Let $M$ be an even-dimensional manifold with a Besse family of hypersurfaces. Then $M$ is closed, and the family is simple Zoll, i.e.\ all its hypersurfaces are regular and embedded. Moreover, if $M$ is orientable, then the hypersurfaces are two-sided $n$-spheres and $M$ is diffeomorphic to $\mathbb{S}^{n+1}$, and otherwise they are one-sided real projective spaces and $M$ is diffeomorphic to $\mathbb{RP}^{n+1}$.
\end{cor}

That $M$ is closed follows from Proposition~\ref{prp:besse_impliese_cmpact}. Moreover, in dimension two the assumptions can be weakened by a theorem of Epstein~\cite{Ep72}, according to which a foliation by circles of a compact $3$-manifold is a Seifert fibration.

\begin{cor}\label{cor:Epstein_refinement}
Let $M$ be a closed surface and let $\mathscr{B}$ be a collection of immersed closed curves in $M$ whose Legendrian lifts to $\PT M$ are the leaves of a foliation. Then $M$ is diffeomorphic to $\Sph^2$ or $\RP^2$, all curves in $\mathscr{B}$ are embedded, and the foliation is given by the fibers of a fiber bundle $\PT M\to\RP^2$.
\end{cor}

In other words, every path geometry on $\mathbb{S}^2$ all of whose paths are closed \cite{BGH95} is simple Zoll. For reversible Finsler metrics on $\mathbb{S}^2$ all of whose geodesics are closed, this also follows from \cite{FLS19, DMMS22}. 

For $n\geq2$, the Legendrian lifts of the hypersurfaces have codimension $n+1\geq3$ in $\PT M$, where compact foliations can have unbounded leaf volume \cite{Su76,EV78}, so that a foliation by compact Legendrian lifts need not be an orbifold fiber bundle a priori. On compact manifolds, Rummler showed that the leaf volume is bounded if the leaves are minimal for some Riemannian metric, and he asked whether compactness is necessary \cite{Ru79}. We prove the following analogue of Wadsley's theorem \cite{Wa75} for minimal hypersurfaces, without assuming $M$ to be compact.

\begin{thm}\label{thm:Wadsley_intro}
Let $(M,g)$ be a connected Riemannian manifold of dimension $n+1\geq3$ and let $\mathscr B$ be a collection of closed, connected, immersed minimal hypersurfaces of $M$ whose Legendrian lifts are embeddings and are the leaves of a smooth foliation of $\PT M$. Then $\mathscr B$ is a Besse family.
\end{thm}

In other words, for minimal hypersurfaces the Besse condition reduces to the analogue of the condition that all geodesics are closed. For closed $M$, Theorem~\ref{thm:Wadsley_intro} follows from Rummler's theorem. For general $M$, it answers Rummler's question for foliations of $\PT M$ by Legendrian lifts of minimal hypersurfaces. Our proof is local. Its first part is analogous to the first half of Wadsley's argument: minimality enters through the first variation formula, which shows that the area of the leaves with trivial holonomy is locally constant, and combined with a counting argument of Epstein \cite{Ep76}, this implies that the linear holonomy of every leaf in the closure of a connected component of the region of locally bounded area is finite. In place of the second half of Wadsley's argument, we use an $\RP^n$-structure that the leaves of such a foliation carry; see Proposition~\ref{prp:proj}. For $n\geq2$, it forces a compact leaf with finite linear holonomy to have finite fundamental group, and hence finite holonomy; see Section~\ref{sec:Wadsley}.

Combining Theorem~\ref{thm:Wadsley_intro} with Theorem~\ref{thm:Berger}, Proposition~\ref{prp:besse_impliese_cmpact} and Corollary~\ref{cor:Gromoll--Grove}, we obtain the following statement.

\begin{cor}\label{cor:Wadsley_all}
Let $(M,g)$ and $\mathscr B$ be as in Theorem~\ref{thm:Wadsley_intro}. If $M$ is closed, then $\mathscr B$ is a Zoll family. If $M$ is even-dimensional, then $M$ is closed, $\mathscr B$ is a simple Zoll family, and $M$ is diffeomorphic to $\Sph^{n+1}$ or $\RP^{n+1}$.
\end{cor}

We record the following open questions.

\begin{qst}
Is there an unoriented transitive family of non-embedded hypersurfaces on an odd-dimensional sphere?
\end{qst}

\begin{qst}\label{qst:top_Wadsley}
Let $M$ be a closed $(n+1)$-manifold with $n\geq2$, and let $\mathcal F$ be a smooth foliation of $\PT M$ by compact Legendrian leaves transverse to the fibers of $\PT M\to M$. Is $\mathcal F$ given by the fibers of an orbifold fiber bundle?
\end{qst}

For $n=1$ the answer is positive by Epstein's theorem \cite{Ep72}. For $n\geq2$, any positive answer must use the contact geometry, since compact foliations in codimension at least three can have unbounded leaf volume \cite{Su76,EV78}. By \cite[Théorème~3]{Ru79}, a positive answer is equivalent to the existence of a Riemannian metric on $\PT M$ for which all leaves are minimal.
\newline

In the following, all manifolds and orbifolds are assumed to be smooth, connected and without boundary, and hypersurfaces are assumed to be immersed, i.e.\ given by immersions $f\colon\Sigma\to M$ of closed manifolds (or orbifolds), unless stated otherwise.

\subsection{Structure of the paper}
In Section~\ref{sec:contact_geometry} we discuss Legendrian fiber bundles and identify Besse families of hypersurfaces as a class of double Legendrian fibrations. In particular, this implies that the regular hypersurfaces of a Besse family are either spheres or real projective spaces. In Section~\ref{sec:Berger_conjec} we make some general observations about Besse families of hypersurfaces and prove Theorem~\ref{thm:Berger}. In Section~\ref{sec:gromoll_grove} we prove the generalized Gromoll--Grove result, Theorem~\ref{thm:Gromoll--Grove}. In Section~\ref{sec:Wadsley} we prove the Wadsley-type theorem, Theorem~\ref{thm:Wadsley_intro}. In Section~\ref{sec:topo_class} we study isotopically reversible submanifolds and provide our proof of the topological classification of manifolds admitting Zoll families of embedded hypersurfaces. Finally, in Section~\ref{sec:deformations} we prove Theorem~\ref{thm:zoll_finite_isometry}, Theorem~\ref{thm:non_embed_odd} and Corollary~\ref{cor:Besse_orbifold} on the existence of nonstandard Zoll resp.\ Besse families of minimal hypersurfaces.
\newline
\newline
\textbf{Acknowledgements.} I thank Marco Radeschi for comments on Lemma~\ref{lem:lemma} and Markus Land for discussions about separating hypersurfaces and spectral sequences.

The proofs in Section~\ref{sec:Berger_conjec} from Lemma~\ref{lem:tau_free} onwards, except Proposition~\ref{Berger_odd_orientable}, until the end of Section~\ref{sec:Wadsley} and of Lemmas~\ref{lem:suf_cond_separating} and~\ref{lem:closed} were developed in interaction with Claude (Anthropic). Claude and ChatGPT (OpenAI) further assisted in literature search, checking drafts and improving the exposition of the paper, including a simplification in Remark~\ref{rem:non-Zoll} suggested by Claude. All mathematical content was checked by the author, who takes full responsibility for it.

\section{Besse families of hypersurfaces via contact geometry}\label{sec:contact_geometry}

In the following we discuss Besse families of hypersurfaces from a contact-geometric perspective. For background on contact geometry, we refer to \cite{Ge08}, and for background on orbifolds to \cite{Th79}, for instance. In particular, recall that the regular part of a connected orbifold, i.e.\ the set of its points with trivial local group, is open, dense and connected \cite[Lemma~5.5]{La18}. For ease of exposition, parts of the  discussion are phrased only in the language of manifolds. Readers primarily interested in the Zoll case may assume all orbifolds to be manifolds. Omitted details in the orbifold case can be filled in by working in local manifold charts and requiring the relevant constructions to be equivariant; cf.\ the proof of Proposition~\ref{prp:contacto_cover}. We will point this out explicitly only where the behavior in the orbifold case may be unexpected and will otherwise use the orbifold versions later without further comment. We frequently use parentheses to indicate the orbifold version of a statement, e.g.\ ``Zoll (Besse)'' and ``(regular) point'' read ``Zoll'' and ``point'' in the manifold case and ``Besse'' and ``regular point'' in the orbifold case. Whether parentheses are used in this sense will be clear from the context.

We begin with a discussion of Legendrian fiber bundles.

\subsection{Legendrian fiber bundles}\label{sub:Lege_fib_bun}

Let $(E^{2n+1},\xi)$ be a contact manifold. A fiber bundle $\mu:E\to M$
is called \emph{Legendrian} if $\ker(T\mu)\subset\xi $ and $\ker(T\mu)$ has rank $n$. This is the maximal possible rank for an integrable distribution contained in the contact distribution. In particular, $M$ has dimension $n+1$. We denote the fibers of $\mu$ by $E_x=\mu^{-1}(x)$ for $x\in M$.

Every manifold $M$ is the base of a Legendrian fiber bundle, namely the projectivized cotangent bundle $\pi_M: \PT M \rightarrow M$ equipped with its natural contact structure $\eta$ defined by $\eta_z=(T_z\pi_M)^{-1}(z)$ for all $z \in \PT M$, where we identify projectivized covectors with hyperplanes. In fact, every Legendrian fiber bundle $\mu :E \rightarrow M$ locally looks like this model \cite[Theorem~2.5]{Pa90}; see also \cite[Appendix~4, p. 367]{Ar89}. More precisely, the projection $\mu$ lifts to a smooth map
\[
  \Phi\colon E\to \PT M,\qquad
  z\mapsto T\mu(\xi_z),
\]
which preserves the contact structures in the sense that $(T\Phi)^{-1}(\eta)=\xi$. Indeed, by definition of $\eta$ and $\Phi$, and since $\ker(T\mu)\subset\xi$, we have
$T\Phi(v)\in\eta_{\Phi(z)} \Leftrightarrow T\mu(v)\in T\mu(\xi_z) \Leftrightarrow v\in\xi_z$
for all $z\in E$ and $v\in T_zE$. Hence $\Phi$ pulls back local contact forms defining $\eta$ to local contact forms
defining $\xi$, and thus volume forms to volume forms, so that it is a local contactomorphism; cf.~\cite[Remark 1.1.6]{Ge08}. Note that in the orbifold setting, such a map need not be locally injective on the underlying spaces, since finite isotropy may still be present.

\subsection{Coorientations and lifts}\label{sub:coor_and_lifts}

A contact manifold $(E^{2n+1},\xi)$ is called \emph{cooriented} if the line bundle $L\coloneqq TE/\xi$ is oriented and \emph{coorientable} if $L$ is orientable. If $(E,\xi)$ is not coorientable, then the orientation double cover $\rho\colon\tilde E\to E$ of $L$ pulls back $\xi$ to a coorientable contact structure $\tilde\xi$ on $\tilde E$. For example, the natural contact structure $\eta$ on the projectivized cotangent bundle $\PT M$ of a manifold $M$ is never coorientable, as the restriction of the corresponding line bundle $L$ to a fiber, isomorphic to the tautological line bundle over a real projective space, is not orientable. The bundle $\PT M$ is doubly covered by the cosphere bundle $\mathbb{S}T^*M\coloneqq (T^* M \setminus 0) /\R_{>0}$ of $M$ on which $\eta$ pulls back to a coorientable contact structure $\tilde \eta$. 

We call a Legendrian fiber bundle $\mu: E \rightarrow M$ \emph{(regularly) coorientable} if the restricted line bundle $L|_{E_x}$ is orientable for one and then every (regular) $x\in M$. In the examples  discussed so far, the Legendrian fiber bundle and its total space are either both coorientable or not. An example of a coorientable Legendrian fiber bundle with non-coorientable total space is given by the Legendrian fibration of $\PT \mathbb{S}^{n+1}$ whose fibers are the Legendrian lifts of the equatorial hyperspheres of $\mathbb{S}^{n+1}$. Its lift to $\ST\mathbb{S}^{n+1}$ generalizes as follows: if $M$ is a Riemannian manifold with injectivity radius $\iota>0$ and $r\in (0,\iota)$, then the Legendrian lifts to $\ST M$ of the distance spheres of radius $r$ in $M$, cooriented by the outward pointing radial vector field, are the fibers of a coorientable Legendrian fiber bundle $\ST M\to M$, whose projection maps each fiber to the center of the corresponding distance sphere; cf.\ \cite[Example~2.5]{Ma26}. Since every manifold admits a metric with positive injectivity radius \cite{Gr78}, such a bundle exists on every manifold. In these examples the fibers are $n$-spheres. For simply connected compact rank one symmetric spaces and suitable radii, these distance spheres are minimal; see \cite{Gu26}.

For a Legendrian fiber bundle $\mu\colon E\to M$, the local contactomorphism $\Phi\colon E\to \PT M$ induces an isomorphism $L\cong\Phi^*L_{\PT M}$ of line bundles. Since $\ST M\to \PT M$ is the orientation double cover of $L_{\PT M}$, the map $\Phi$ lifts to a local contactomorphism to $\ST M$ if and only if $(E,\xi)$ is coorientable. In this case we also write $\tilde E=E$ and $\rho=\mathrm{id}$. In any case, $\Phi\circ\rho\colon\tilde E\to \PT M$ lifts to a local contactomorphism $\tilde\Phi\colon\tilde E\to\ST M$. Likewise, for a (regular) point $x\in M$ the restriction $\Phi_x\colon E_x\to \PT _xM$ lifts to $\ST _xM$ if and only if the Legendrian fiber bundle $\mu$ is (regularly) coorientable.

If the fibers $E_x$ are closed, then the local diffeomorphisms $\Phi_x$ are covering maps. If the fibers are in addition connected and $n\geq 1$, then this implies that the (regular) fibers $E_x$ are diffeomorphic either to the $n$-sphere or to the $n$-dimensional real projective space depending on whether $\mu$ is (regularly) coorientable or not. Since the degree of the covering $\Phi_x$ is independent of a (regular) point $x$, in this case also the local diffeomorphism $\Phi: E \rightarrow \PT M$ is a covering map of the same degree. In particular, if $n\geq 2$ and $\mu$ is not (regularly) coorientable, then the degree is one and $\Phi: E \rightarrow \PT M$ is a contactomorphism. If $n\geq 2$ and $(E,\xi)$ is coorientable, then  $\tilde \Phi : E \rightarrow \ST M$ is a contactomorphism. 

In the following proposition we record the orbifold version of these observations and provide further details on the proof in this setting; see \cite[Theorem~6.2 and Corollary~6.3]{Pa90} for related statements on Legendre foliations and \cite[Propositions~3.1 and~3.4 and Corollary~3.2]{Ma26} for related statements for unoriented transitive families of closed embedded hypersurfaces.

\begin{prp}\label{prp:contacto_cover} Let $(E,\xi)$ be a contact manifold and $\mu\colon E\to\Orb$ a Legendrian orbifold fiber bundle with closed connected fibers. Then $\Phi\colon E\to\PT\Orb$ is a covering map. Moreover, if the base has dimension at least three, then the following hold.
\begin{enumerate}
\item If $\mu$ is not regularly coorientable, then $\Phi$ is a contactomorphism.
\item If $(E,\xi)$ is coorientable, then $\Phi$ lifts to a contactomorphism $\tilde\Phi\colon E\to\ST\Orb$.
\end{enumerate}
\end{prp}
\begin{proof} The first claim is local, so we can assume that $\Orb=M/\Gamma$ where $\Gamma$ is a finite group that acts smoothly on a manifold $M$. The bundle $\mu$ is then covered by a $\Gamma$-equivariant bundle $ \hat \mu: \hat E\rightarrow M$, and the bundle $\PT \Orb \rightarrow \Orb$ is covered by the $\Gamma$-equivariant bundle $\PT M\rightarrow M$. The argument from the manifold case gives a $\Gamma$-equivariant covering map $\hat \Phi: \hat E \rightarrow \PT M$, which descends to an orbifold covering map $\Phi: E \rightarrow \PT \Orb$.

If $\mu$ is not regularly coorientable, then the covering map $\Phi$ restricts to a diffeomorphism between fibers over regular points in $\Orb$. Since an orbifold covering map is a diffeomorphism if some regular point in the target has only one preimage, claim (i) follows.

In case (ii) we choose in an equivariant model as above the lift of $\hat\Phi$ that maps the coorientation of $\hat\xi$ pulled back from $E$ to that of the natural contact structure on $\ST M$. Since $\Gamma$ preserves these coorientations, the lift is $\Gamma$-equivariant and descends to a covering $\tilde\Phi\colon E\to\ST\Orb$. Since $n\geq 2$, it restricts to a diffeomorphism over a regular point in $\Orb$, and hence $\tilde\Phi$ is a diffeomorphism as before.
\end{proof}

The remaining case in which $\mu$ is regularly coorientable but $(E,\xi)$ is not coorientable will not be needed in the following.

We recall another standard contact geometry fact; cf.~\cite[Chapter~5, \S 1, Remark~1]{AG01}.

\begin{lem}\label{lem:contact_equiv}
Suppose that, in the setting of Proposition~\ref{prp:contacto_cover}, the contact manifold $(E,\xi)$ is coorientable, so that $\Phi$ lifts to $\tilde\Phi\colon E\to\ST\Orb$, and let $F$ be a contactomorphism of $(E,\xi)$ that preserves the fibers of $\mu$. Then $F$ descends to a diffeomorphism $f$ of $\Orb$ with $\tilde\Phi\circ F=Tf\circ\tilde\Phi$ or $\tilde\Phi\circ F=-Tf\circ\tilde\Phi$, depending on whether $F$ preserves the coorientation of $\xi$ or not.
\end{lem}

\begin{proof} That the induced map $f$ is a diffeomorphism follows from the local structure of $\mu$.

By construction, for $z \in E$, $\tilde \Phi(z)$ is the hyperplane $T\mu(\xi_z)$ cooriented by the coorientation of $\xi_z$ transferred to $T\mu(\xi_z)$ via $T\mu$. The claim follows from $\mu\circ F=f\circ \mu$ and the fact that both $\tilde \Phi$ and $Tf$ preserve coorientations.
\end{proof}

\subsection{Zoll and Besse Legendrian fiber bundles}\label{sub:double_Legendrian}

We now give a contact-geometric interpretation of the notion of a Besse family of hypersurfaces; cf.\ \cite{BGH95,LM02,La20} for double fibrations in the geodesic case. We call a Legendrian fiber bundle $\mu:E \rightarrow M$
of a contact manifold $(E,\xi)$ a \emph{Zoll (Besse) Legendrian fiber bundle} if all its fibers are connected and closed, and if there exists a Legendrian (orbifold) fiber bundle $\nu\colon E\to N$ whose fibers are connected, closed and transverse to the fibers of $\mu$ in $\xi$, i.e.\ $\ker T_z\mu\cap\ker T_z\nu=0$, or equivalently $\ker T_z\mu\oplus\ker T_z\nu=\xi_z$, for all $z\in E$. In this case the pair $(\mu,\nu)$ defines a \emph{double Legendrian fibration}. Note that the definition becomes symmetric in $\mu$ and $\nu$ if we also allow $M$ to be an orbifold. With this convention, $\nu$ is a Besse Legendrian fiber bundle, and it is Zoll if $M$ is a manifold.

Given such a double Legendrian fibration, the map $(\mu,\nu)\colon E\to M\times N$ is an immersion, as $\ker T\mu\cap\ker T\nu=0$, and hence locally injective. For the same reason, the fibers of $\mu$ immerse into $N$ and the fibers of $\nu$ immerse into $M$. Moreover, if $M$ and $N$ are manifolds, then all fibers of $\mu$ embed into $N$ if and only if all fibers of $\nu$ embed into $M$. If this holds, we call the double Legendrian fibration and the corresponding families of hypersurfaces \emph{simple Zoll}.

We denote by $\Sigma_y$ the hypersurface $\mu|_{E_y}\colon E_y\to M$, $y\in N$, and analogously by $\Sigma_x$ the hypersurface $\nu|_{E_x}\colon E_x\to N$, $x\in M$. A hypersurface in a manifold is called \emph{two-sided} if its normal bundle is orientable. For a (regular) point $y\in N$, the hypersurface $\Sigma_y$ in $M$ is two-sided if and only if $\nu$ is (regularly) coorientable.

For $n \geq 2$, an unoriented (Besse) transitive family of hypersurfaces, i.e.\ a Zoll (Besse) family, on an $(n+1)$-manifold $M$ is the same as a Zoll (Besse) Legendrian fiber bundle $\mu:E \rightarrow M$ that is not (regularly) coorientable. Indeed, by Proposition~\ref{prp:contacto_cover}, (i), we can identify $\mu:E \rightarrow M$ with the natural bundle $\PT M \rightarrow M$ in this case. Moreover, the (regular) hypersurfaces in the family are two-sided if and only if the corresponding Legendrian fiber bundle $\nu$ is (regularly) coorientable. We record the following statement for later use.

\begin{lem}\label{lem:simple_charact}
A Zoll family $\mathscr{Z}$ on an $(n+1)$-manifold $M$, $n\geq 1$, is simple Zoll if one of its hypersurfaces is embedded.
\end{lem}
\begin{proof} We use an open-closed argument as in \cite{GG81}. A hypersurface in $\mathscr{Z}$ can only have transverse self-intersections, since its  Legendrian lift to  $\PT M$ is embedded. Since transverse intersections persist under $C^1$-small deformations, the set of non-embedded hypersurfaces in $P$ is open. This set is also closed. Indeed, suppose a sequence of non-embedded hypersurfaces $\Sigma_k$ in $\mathscr{Z}$ converges in $P$ to a hypersurface $\Sigma$ in $\mathscr{Z}$. We can then find points $z_k\neq w_k$ in the Legendrian lift of $\Sigma_k$ with $\mu(z_k) = \mu(w_k)$. By local triviality of $\nu$ and compactness of the hypersurfaces, the preimage $\nu^{-1}(K)$ of a compact neighborhood $K$ of $\Sigma$ in $P$ is compact.  Therefore, after passing to a subsequence, we can assume that $z_k$ and $w_k$ converge to points $z$ resp.\ $w$ in the Legendrian lift of $\Sigma$ with $\mu(z)=\mu(w)$. By local injectivity of $(\mu,\nu)$ we have $z\neq w$ and so $\Sigma$ is not embedded. Since $P$ is connected, the claim follows.
\end{proof}

\section{Berger-type conjecture for Besse families}\label{sec:Berger_conjec}

Let $n\geq 1$. For a Besse family of hypersurfaces on an $(n+1)$-manifold $M$, we consider the double cover $\ST M\to\PT M$ with deck transformation $\tau$. For a regular point $\sigma\in P$, we identify the hypersurface $\Sigma=\Sigma_\sigma$ with its Legendrian lift $E_\sigma=\nu^{-1}(\sigma)\subset\PT M$ and denote the connected components of its preimage in $\ST M$ by $\Sigma^+$. Together with the connected components of the preimages of the exceptional fibers, the $\Sigma^+$ form the fibers of a Legendrian orbifold fiber bundle $\nu^+\colon\ST M\to P^+$ with connected fibers, as can be seen by pulling back the double cover $\ST M\to\PT M$ to local models of $\nu$; cf.\ the beginning of Section~\ref{sec:contact_geometry}. We obtain a commutative diagram 
\[
\begin{tikzcd}
\Sigma^+ \arrow[d] \arrow[r]  & \ST M \arrow[d] \arrow[r, "\nu^+"]  &  P^+  \arrow[d] \\
 \Sigma \arrow[r]  & \PT M \arrow[r, "\nu"]   &   P 
\end{tikzcd}
\]
in which the horizontal sequences are orbifold fiber bundles and the vertical maps are orbifold coverings. By the discussion in Section~\ref{sec:contact_geometry} we see that if the regular hypersurfaces are two-sided, i.e.\ if $\nu$ is regularly coorientable, then the left vertical covering is a diffeomorphism between $n$-spheres and the right vertical covering is two-fold. Otherwise, the left vertical covering is a two-fold covering from an $n$-sphere to an $n$-dimensional real projective space and the right vertical covering is a diffeomorphism. The deck transformation $\tau$ descends to an involution $\tau^+$ of $P^+ $ with quotient $P$.

The same argument applies to arbitrary finite coverings: if $M'\to M$ is a finite covering, then the connected components of the preimages of the fibers of $\nu$ under the induced covering $\PT M'\to\PT M$ are the fibers of a Legendrian orbifold fiber bundle $\nu'\colon\PT M'\to P'$, as can be seen by pulling back the covering to local models of $\nu$. Hence every Besse family on $M$ lifts to a Besse family on $M'$.

We call an orbifold fiber bundle orientable if its fibers can be endowed with consistent orientations.

\begin{lem}\label{lem:lifted_orientations}
For $n\geq 1$  consider a Besse family of hypersurfaces on an orientable $(n+1)$-manifold $M$. Then also $P^+ $ and $\nu^+$ are orientable. Moreover, $\tau$ reverses the orientation of the fibers of $\nu^+$ and it preserves the orientation of $P^+ $ if and only if the dimension $n+1$ of $M$ is odd. In particular, if $n+1$ is even, then $\nu$ is regularly coorientable and its regular fibers are spheres.
\end{lem}

\begin{proof}
We can use the coorientation of the natural contact structure and an orientation of $M$ to choose consistent orientations for the fibers of $\nu^+$. Since also the total space $\ST M$ is orientable, we can then define an orientation on $P^+$.

Since $\tau$ covers the identity of $M$ and flips the coorientation of the contact structure, it reverses the orientations of the fibers of $\nu^+$. Moreover, $\tau$ acts as the antipodal map on the fibers of $\ST M\to M$, so it preserves the orientation of $\ST M$ if and only if $n+1$ is even. Hence $\tau^+$ preserves the orientation of $P^+$ if and only if $n+1$ is odd. In particular, if $n+1$ is even, then $\tau^+$ acts nontrivially on $P^+$, so $\nu$ is regularly coorientable and its regular fibers are spheres, as discussed in Subsection~\ref{sub:coor_and_lifts}.
\end{proof}

The following lemma will in particular allow us to deal with the odd-dimensional case. 

\begin{lem}\label{lem:orbifold_manifold}
Let $\Orb$ be an orbifold. If $\ST \Orb$ or $\PT \Orb$ is a manifold, then all singular points of $\Orb$ are isolated. If $\Orb$ is in addition odd-dimensional, then the only possible nontrivial local group is $\{\pm \mathrm{id}\}$. If $\Orb$ is in addition odd-dimensional and orientable, then $\Orb$ is a manifold.
\end{lem}
\begin{proof}  
If $\PT \Orb$ is a manifold, then so is its cover $\ST \Orb$. The latter is a manifold if and only if every local group acts freely on the unit sphere of its tangent space representation, which in turn is equivalent to the statement that all singular points of $\Orb$ are isolated. The second claim follows since the only nontrivial free linear action on an even-dimensional sphere is the antipodal action. The last claim follows since the antipodal action of an even-dimensional sphere reverses the orientation.
\end{proof}

The next statement will be applied several times.

\begin{lem}\label{lem:fix}
Let $N$ be a connected manifold, $f$ a smooth involution of $N$ and
$\varepsilon\in\{\pm 1\}$. Suppose that $\varepsilon\, Tf$ acts freely on $\ST N$.
\begin{enumerate}[label=(\alph*)]
\item If $\varepsilon=-1$ and $f$ has a fixed point, then $f=\mathrm{id}$.
\item If $\varepsilon=+1$, then $T_xf=-\mathrm{id}$ at every fixed point $x$ of $f$. In
particular, the fixed points are isolated.
\end{enumerate}
\end{lem}
\begin{proof}
At a fixed point $x$ the differential $d_xf$ is a linear involution, and freeness
means that it has no eigenvalue $\varepsilon$. Hence $T_xf=-\varepsilon\,\mathrm{id}$. If
$-\varepsilon=1$, then $f$ is an isometry of an $f$-invariant Riemannian metric
fixing $x$ with trivial differential, so $f=\mathrm{id}$ since $N$ is connected.
\end{proof}

\begin{lem}\label{lem:tau_free}
For $n\geq 1$ consider a Besse family of hypersurfaces on an orientable $(n+1)$-manifold $M$ with $\nu$ regularly coorientable. Then $\tau^+$ acts freely on $P^+$.
\end{lem}
\begin{proof} Suppose that $\tau^+$ fixes a point $\sigma \in P^+$.

We first assume that $n\geq 2$. In this case, Proposition~\ref{prp:contacto_cover}, (ii), applied to $\nu^+$ gives a contactomorphism $\tilde \Phi^+\colon \ST M\to \ST P^+$ with $\tilde \Phi^+\circ \tau = - T\tau^+\circ \tilde \Phi^+$
by Lemma~\ref{lem:contact_equiv}. In particular, $- T\tau^+$ acts freely on $\ST  P^+$ as $\tau$ acts freely on $\ST M$. Then no lift of $T_\sigma \tau^+$ to $\mathrm{O}(n+1)$ has an eigenvalue $-1$. For $n+1$ even this is not possible, since every such lift has determinant $-1$ by Lemma~\ref{lem:lifted_orientations}, and an orthogonal transformation  with determinant $-1$ has eigenvalue $-1$. For $n+1$ odd, $P^+$ is a manifold by Lemma~\ref{lem:orbifold_manifold} since $\ST P^+ \cong \ST M$ is a manifold and $P^+$ is orientable by Lemma~\ref{lem:lifted_orientations}. Lemma~\ref{lem:fix} (a) applied to $f=\tau^+$  implies that $\tau^+=\mathrm{id}$. But $\tau^+$ acts transitively on the fibers of the covering $P^+ \rightarrow P$, which is two-fold as $\nu$ is regularly coorientable, a contradiction.

Now let $n=1$. Then $\tau$ restricts to a free involution of the fiber of $\nu^+$ over $\sigma$, which is an embedded circle in $\ST M$. As in the proof of Lemma~\ref{lem:lifted_orientations}, we see that $\tau$ reverses its orientation. Since $\Sph^1$ does not admit an orientation-reversing free involution, we again obtain a contradiction.
\end{proof}

The following observation will be used to show compactness and the rational homology sphere property in even dimensions; see Propositions~\ref{prp:besse_impliese_cmpact} and~\ref{prp:even_rat_hom_sphere}, respectively.

\begin{lem}\label{lem:induced_involution_cohomology}
Let $M$ be an even-dimensional manifold. Then the fiberwise involution $\tau$ of $E=\ST M$ acts trivially on $H^*(E;\Q)$ and $H_*(E;\Q)$.
\end{lem}
\begin{proof}
Since $\tau$ covers the identity of $M$ and preserves the orientation of the fibers of $\mu^+$, it acts trivially on the second, and then also the infinity page of the cohomological Serre spectral sequence of $\mu^+$. Thus $\tau^*$ preserves the finite Serre filtration of $H^k(E;\Q)$ and acts trivially on its associated graded pieces, so $(\tau^*-1)^{k+1}=0$. Since also $(\tau^*)^2=1$ and $\gcd((x-1)^{k+1},x^2-1)=x-1$ over $\Q$, we get $\tau^*=\mathrm{id}$. The same argument applies in homology.
\end{proof}

In the following, we use the Euler class $e(p)\in H^{n+1}(B;\Z)$ and the Gysin sequence of an orientable $n$-sphere bundle $p\colon E\to B$; see \cite[\S4.D]{Hat02}, for instance. For an orbifold sphere bundle, we apply this to the Borel construction. Recall that a model of the classifying space of an orbifold $\Orb$ of dimension $n+1$ is $B\Orb=\mathrm{Fr}(\Orb)\times_{\mathrm{O}(n+1)}E\mathrm{O}(n+1)$, where $\mathrm{Fr}(\Orb)$ denotes the orthonormal frame bundle (with respect to some Riemannian metric); see e.g.\ \cite[Section~2.5]{ALR21} and \cite[Section~3.1]{KL21}. An orbifold $n$-sphere bundle $p\colon E\to\Orb$ whose total space $E$ is a manifold then induces an $n$-sphere bundle over $B\Orb$ with total space homotopy equivalent to $E$; see \cite[Section~3.1]{KL21}, whose discussion on circle bundles applies here analogously. If the fibers of $p$ are consistently oriented, then these orientations pull back to an $O(n+1)$-invariant orientation of the fibers of the pullback of $p$ to the frame bundle of $\Orb$, which therefore descends to the Borel construction. Hence the induced fibration is orientable and we write $e(p)\in H^{n+1}(B\Orb;\Z)$ for its Euler class. Moreover, the natural map $B\Orb\to|\Orb|$ induces an isomorphism in rational cohomology \cite[Proposition~2.12]{ALR07}.

\begin{prp}\label{prp:besse_impliese_cmpact}
Let $M$ be a manifold of even dimension $n+1\geq 2$, which  carries a Besse family of hypersurfaces. Then $M$ and $P^+$ are compact.
\end{prp}
\begin{proof}
By passing to the orientable double cover and lifting the Besse family we can assume that $M$ is orientable. Since both the fibers of $\mu^+$ and $\nu^+$ are compact, $M$ is compact if and only if $P^+$ is compact.

Suppose that $M$ and $P^+$ are non-compact, let $E=\ST M$ and let $[\Sigma] \in H_n(E;\Q)$ be the class of the consistently oriented fibers of $\nu^+$; see Lemma~\ref{lem:lifted_orientations}. Since $\tau$ reverses the orientations of the fibers by Lemma~\ref{lem:lifted_orientations}, we have $\tau_* [\Sigma]= - [\Sigma]$. Since $\tau$ acts trivially on $H_*(E;\Q)$ by Lemma~\ref{lem:induced_involution_cohomology}, it follows that $[\Sigma]=0$.

We have $H^{n+1}(BP^+;\Q)\cong H^{n+1}(|P^+|;\Q)=0$, since $|P^+|$ is a non-compact connected rational homology manifold of dimension $n+1$ (the local groups act orientation-preservingly, since $P^+$ is orientable by Lemma~\ref{lem:lifted_orientations}), and the top-dimensional rational homology, and hence cohomology, of such a space vanishes; cf.\ \cite[Proposition~3.29]{Hat02}. Therefore, the Euler class of $\nu^+$ vanishes. Since $e(\nu^+)=0$, the map $H^n(E;\Q)\to H^0(BP^+;\Q)\cong\Q$ in the Gysin sequence of $\nu^+$, which is given by evaluation on the class of a regular fiber $\Sigma^+$ of $\nu^+$, is surjective. Hence some class $\varphi\in H^n(E;\Q)$ satisfies $\varphi([\Sigma])\neq0$, a contradiction.
\end{proof}

\subsection{Orientable manifolds of odd dimension}\label{sub:Berger_orientable_odd}

We obtain the following statement in odd dimensions.

\begin{prp} \label{Berger_odd_orientable}
Let $M$ be an orientable odd-dimensional manifold. Then every Besse family of hypersurfaces on $M$ is Zoll, i.e. all its hypersurfaces are regular.
\end{prp}
\begin{proof} By Lemma~\ref{lem:lifted_orientations}, $P^+$ is orientable and $\tau^+$ preserves its orientation, so $P$ is orientable as well. Since $M$ is odd-dimensional and $n\geq 1$, its dimension is at least three, so Proposition~\ref{prp:contacto_cover} applies. If $\nu$ is not regularly coorientable, then $\Phi: \PT M\rightarrow \PT P$ is a diffeomorphism by Proposition~\ref{prp:contacto_cover}, (i), and so $P$ is a manifold by Lemma~\ref{lem:orbifold_manifold}.

If $\nu$ is regularly coorientable, then  $\tilde \Phi^+: \ST M \rightarrow \ST P^+$ is a diffeomorphism by Proposition~\ref{prp:contacto_cover}, (ii), applied to $\nu^+:\ST M \rightarrow P^+$, whose total space is coorientable, and so $P^+$ is a manifold by Lemma~\ref{lem:orbifold_manifold}. Since the action of $\tau^+$ on $P^+$ is free by Lemma~\ref{lem:tau_free} also $P=P^+ /\tau^+$ is a manifold.
\end{proof} 

\subsection{Orientable manifolds of even dimension}\label{sub:Berger_orientable_even}

We continue to work with the Gysin sequence and the Euler class of an orbifold sphere bundle as discussed before Proposition~\ref{prp:besse_impliese_cmpact}.

\begin{prp}\label{prp:even_rat_hom_sphere}
Let $M$ be an orientable manifold of even dimension $n+1 \geq 2$ with a Besse family of hypersurfaces. Then both $M$ and the underlying space $|P^+|$ are rational homology spheres. In particular, $H^1(M;\Z)=0$ and $H^1_{orb}(P^+;\Z)=0$.
\end{prp}
\begin{proof}
By Proposition~\ref{prp:besse_impliese_cmpact} $M$ and $P^+$ are compact. By Lemma~\ref{lem:lifted_orientations} the orbifold fiber bundle $\nu^+:E \rightarrow P^+$ is orientable. The rational Gysin sequence of the associated Borel construction shows that the $\tau$-equivariant map
\[
        (\nu^+)^*: H^j(BP^+;\Q) \rightarrow H^j(E;\Q)
\]
is an isomorphism for $j\leq n-1$ and injective for $j=n$. By Lemma~\ref{lem:induced_involution_cohomology}, $\tau^+$ hence also acts trivially on $H^j(BP^+;\Q)$ for all $j\leq n$. The natural projection $BP^+\to|P^+|$ is $\tau^+$-equivariant and induces isomorphisms in rational cohomology \cite[Proposition~2.12]{ALR07}, so $\tau^+$ also acts trivially on $H^j(|P^+|;\Q)$ for $j\leq n$. On the other hand, $(\tau^+)^*=-\mathrm{id}$ on $H^{n+1}(|P^+|;\Q)\cong\Q$, since $\tau^+$ reverses the orientation of $P^+$ by Lemma~\ref{lem:lifted_orientations}.

For $1\leq j \leq n$ both $j$ and $n+1-j$ are at most $n$, so for $x\in H^j(|P^+|;\Q)$ and $y \in H^{n+1-j}(|P^+|;\Q)$ we have
\[
    x \cup y =\tau^*(x \cup y) = - x \cup y
\]
and hence $x \cup y=0$. Since $|P^+|$ is a closed oriented rational homology manifold \cite{Sa56}
it satisfies Poincar\'e duality \cite[Chapter~II, \S 2]{Bo60}
and so $H^j(|P^+|;\Q)=0$ for $1 \leq j\leq n$. The Gysin sequence of $\mu^+$ gives $H^j(M;\Q) \cong H^j(E;\Q)$ for $j\leq n-1$. If $n+1 \geq 4$, then this together with Poincar\'e duality on $M$ implies that $H^j(M;\Q)=0$ for $1 \leq j \leq n$.

Assume now that $n+1=2$. The Gysin sequence of $\nu^+$ gives
\[
0 = H^1(|P^+|; \Q) \rightarrow H^1(E; \Q) \rightarrow H^0(|P^+|; \Q) \cong \Q
\]
so that $b_1(E) \leq 1$. The same portion of the Gysin sequence for $\mu^+$ shows that $b_1(E) \geq b_1(M)$. Since the first Betti number of an orientable surface is even, it follows that $M$ is topologically a $2$-sphere. 

Since $H_1(M;\Q)=0$ and $H_1^{orb}(P^+;\Q)=0$, the groups $H_1(M;\Z)$ and $H_1^{orb}(P^+;\Z)$ are torsion, so $H^1(M;\Z)\cong\mathrm{Hom}_\Z(H_1(M;\Z),\Z)=0$ and likewise $H^1_{orb}(P^+;\Z)=0$.
\end{proof}

We obtain the following consequences.

\begin{cor}\label{cor:euler_class}
Let $M$ be an orientable manifold of even dimension $n+1 \geq 2$ with a Besse family of hypersurfaces. Then 
\[
\Z_2\cong H^{n+1}(M;\Z)/\langle e(\mu^+) \rangle \cong H^{n+1}(E;\Z) \cong H^{n+1}(BP^+;\Z)/\langle e(\nu^+) \rangle .
\]
\end{cor}

\begin{proof}
The Gysin sequence of $\mu^+$,
\begin{equation}\label{eq:gysin}
  H^0(M;\Z) \xrightarrow{\ \cup\, e\ } H^{n+1}(M;\Z)
    \xrightarrow{\ (\mu^+)^*\ } H^{n+1}(E;\Z) \longrightarrow H^1(M;\Z) = 0,
\end{equation}
identifies $H^{n+1}(E;\Z)$ with $\ H^{n+1}(M;\Z)/\langle e(\mu^+) \rangle$. Since the Euler class of the cosphere bundle $\mu^+$, which is isomorphic to the unit tangent bundle via any Riemannian metric, evaluates to $\chi(M) = 2$ on a fundamental class of the even-dimensional rational homology sphere $M$ (see \cite[Corollary~11.12]{MS74}), we have $H^{n+1}(E;\Z) \cong \Z_2$. The same computation for $\nu^+$ completes the proof.
\end{proof}

\begin{cor}\label{cor:no_finite_covers}
Let $M$ be an orientable manifold of even dimension $n+1 \geq 2$ with a Besse family of hypersurfaces. Then $\pi_1(M)$ and $\pi_1(|P^+|)$ have no proper finite index subgroups.
\end{cor}
\begin{proof}
Suppose $M'\rightarrow M$ is a finite covering of degree $d$. Then the Besse family on $M$ lifts to a Besse family on $M'$. By Proposition~\ref{prp:even_rat_hom_sphere} both $M$ and $M'$ are rational homology spheres. It follows that $2=\chi(M')=d\chi(M)=2d$ and hence $d=1$.

If $n+1=2$, then $|P^+|$ is simply connected by Proposition~\ref{prp:even_rat_hom_sphere}. If $n+1 \geq 4$, then $\pi_1(M) \cong \pi_1^{\mathrm{orb}}(P^+)$ by the long exact homotopy sequences of $\mu^+$ and $\nu^+$. Since $\pi_1(|P^+|)$ is a quotient of $\pi_1^{\mathrm{orb}}(P^+)$, the claims follow.
\end{proof}

\begin{prp}\label{prp:berger_even_dim_orient}
Let $M$ be an orientable manifold of even dimension $n+1 \geq 2$ with a Besse family of hypersurfaces. Then the Besse family is Zoll.
\end{prp}
\begin{proof} By Lemma~\ref{lem:tau_free} it suffices to show that $P^+$ is a manifold. If $P^+$ has no singular points, there is nothing to show.

We first assume that $n+1 \geq 4$. By Proposition~\ref{prp:contacto_cover} applied to $\nu^+$ and Lemma~\ref{lem:orbifold_manifold} all singular points of $P^+$ are isolated. By compactness of $P^+$, there are only finitely many, which we denote by $\sigma_1,\ldots, \sigma_k$. Let $\Gamma_i$, $i=1,\ldots,k$, denote the corresponding local groups and let $U_i \cong D^{n+1}/\Gamma_i$ be disjoint neighborhoods of these singular points. We set $P_0:= P^+ \backslash \bigcup_i U_i^{\circ}$. Then $BP^+$ is the union of $BP_0 \simeq P_0$ and the spaces $BU_i\simeq B\Gamma_i$ glued along spaces homotopy equivalent to the closed $n$-manifolds $\partial U_i=\Sph^n/\Gamma_i$. Since $\Gamma_i$ acts freely and orientation-preservingly  on $\Sph^n$ it has periodic cohomology of period dividing $n+1$, and so $H^n(B\Gamma_i;\Z)=0$ and $H^{n+1}(B\Gamma_i;\Z)\cong \Z_{|\Gamma_i|}$ \cite[Chapter VI, \S4 and \S 9]{Br94}. Moreover, $H^{n+1}(P_0;\Z)=0$, and, by the long exact sequence of the pair $(P_0,\partial P_0)$, the image of the inclusion induced map $H^n(P_0;\Z) \rightarrow H^n(\partial P_0;\Z)$ is the kernel of the sum map $\Z^k\cong H^n(\partial P_0;\Z)\rightarrow H^{n+1}(P_0,\partial P_0;\Z)$. Therefore, Mayer--Vietoris gives a short exact sequence
\[
0\to\Z\, \alpha\xrightarrow{\ \delta^*\ } H^{n+1}(BP^+;\Z)\xrightarrow{\ q\ }\bigoplus_{i=1}^k\Z_{|\Gamma_i|}\to 0 ,
\]
where $\alpha \in H^n(\partial P_0;\Z)$ is a lift of a generator $\beta \in H^{n+1}(P_0,\partial P_0;\Z)$. 

By excision and contractibility of the cones $|U_i|$, we have isomorphisms
\begin{equation}\label{eq:exci_cont}H_{n+1}(P_0,\partial P_0;\Z)\cong H_{n+1}(|P^+|;\Z), \quad H^{n+1}(P_0,\partial P_0;\Z)\cong H^{n+1}(|P^+|;\Z).
\end{equation}
The image $[|P^+|]\in H_{n+1}(|P^+|;\Z)$ of a relative fundamental class $[P_0] \in H_{n+1}(P_0,\partial P_0;\Z)$ under the first isomorphism in \eqref{eq:exci_cont} is a fundamental class in the sense that it restricts to a generator of the local homology $H_{n+1}(|P^+|,|P^+|\setminus\{x\};\Z)\cong\Z$ at every regular point $x \in P^+$ and a nonzero element in $H_{n+1}(|P^+|,|P^+|\setminus\{x\};\Q)\cong\Q$ at every point $x\in P^+$.

Since the projection $BP^+\to|P^+|$ is compatible with the decompositions into $P_0$ and the contractible cones $|U_i|$, naturality of the Mayer--Vietoris sequences shows that the induced map $H^{n+1}(|P^+|;\Z)\rightarrow H^{n+1}(BP^+;\Z)$
sends the image of $\alpha$ under the Mayer--Vietoris coboundary of the decomposition of $|P^+|$ to $\delta^*\alpha$. This image coincides  up to sign with the image of $\beta$ under the second isomorphism in \eqref{eq:exci_cont}. Hence, evaluating classes in $H^{n+1}(BP^+;\Q)$ on $[|P^+|]$ via the isomorphism $H_{n+1}(BP^+;\Q)\cong H_{n+1}(|P^+|;\Q)$ induced by the projection,  we have $\langle\delta^*\alpha,[|P^+|]\rangle=\pm\langle\beta,[P_0]\rangle=\pm1$.

Integrally, we have $H^{n+1}(BP^+;\Z)\cong \Z \oplus T$ for some torsion group $T$. Let $u \in H^{n+1}(BP^+;\Z)$ be a lift of a generator of $H^{n+1}(BP^+;\Z)/T$ and let $l,m \in \Z \backslash \{0\}$ with $e(\nu^+)= l u$ and $\delta^*(\alpha)=mu$ rationally.

Since $\nu^+$ is isomorphic to the cosphere bundle of $P^+$ by Proposition~\ref{prp:contacto_cover}, (ii), its Euler class is the Euler class of the orbifold tangent bundle of $P^+$. With real coefficients, the equivariant de Rham theorem for the locally free action of $\mathrm{O}(n+1)$ on the frame bundle of $P^+$ identifies Borel cohomology with orbifold de Rham cohomology \cite[Theorem~2.5.1 and Theorem~4.3.1]{GS99}, and by Chern--Weil theory the Euler class corresponds to the class of the Euler form under this identification (cf.\ the generalized Gauss--Bonnet theorem in \cite[Appendix~C]{MS74} for vector bundles over manifolds, and \cite[Section~3.2 and Lemma~3.5]{KL21} for circle bundles over orbifolds). Hence $\chi_{\mathrm{orb}}(P^+)=\langle e(\nu^+),[|P^+|]\rangle=\pm l/m$ by \cite{Sa57}.

Moreover, we have $l|T|=2$ by Corollary~\ref{cor:euler_class} and $m|T|= \prod_i |\Gamma_i|$ by the above short exact sequence. Together with $\chi_{\mathrm{orb}}(P^+) = \chi(|P^+|) - \sum_i \left(1- \frac{1}{|\Gamma_i|} \right)$ \cite{Sa57} and $\chi(|P^+|)=2$ we obtain
\[
    \frac{\pm2}{\prod_i |\Gamma_i|} =  \chi_{\mathrm{orb}}(P^+) = 2 -  \sum_i \left( 1 - \frac{1}{|\Gamma_i|} \right)  .
\]
By Lemma~\ref{lem:tau_free} $\tau^+$ permutes singular points without fixed points while preserving the orders of the local groups. Writing $d_1,d_1,\ldots,d_r,d_r$ for these orders with $d_i \geq 2$ and $k=2r$, we obtain
\[
    \pm1  = \prod_i d_i \left((1-r)  \prod_i d_i + \sum_i  \prod_{j\neq i} d_j \right)   .
\]
Since the term in the bracket is an integer, we obtain a contradiction. Hence, $P^+$ is a manifold. 

Let now $n+1=2$. In this case $M$ is a $2$-sphere by Proposition~\ref{prp:even_rat_hom_sphere}. Hence, $\PT M \cong L(4,1)$  and $\ST M\cong\mathbb{RP}^3$ are lens spaces and the fiber bundles define Seifert fibrations on them. Since $P^+ $ is orientable by Lemma~\ref{lem:lifted_orientations}, the classification of Seifert fibrations on $\mathbb{RP}^3$ implies that $P^+ $ is a spindle orbifold $S^2(p,q)$ with at most two singular points of order $p$ and $q$ satisfying $\gcd (p,q) \in \{1,2\}$; see \cite[Theorem~5.4]{FLS19}, for instance. Since $\tau^+$ acts freely on $P^+$ by Lemma~\ref{lem:tau_free}, we see that singular points can only occur as a pair and that they are interchanged by $\tau^+$. This implies $p=q \in \{1,2\}$. By \cite[Lemma~3.6]{La20} there is no Seifert fibration on $\mathbb{RP}^3$ with orbit space $S^2(2,2)$; cf.\ \cite[Example~4.18]{GL18}. Hence, we see that $P^+\cong \Sph^2$ and $P\cong \mathbb{RP}^2$ are manifolds.
\end{proof}

\subsection{Non-orientable case}\label{sub:Berger_non_orientable}

Let $M$ be a non-orientable closed manifold and suppose that it carries a Besse family. The family lifts to a Besse family on the orientable double cover $\tilde M$ of $M$ with parameter spaces $\tilde P^+$ and $\tilde P$. By Propositions~\ref{prp:berger_even_dim_orient} and~\ref{Berger_odd_orientable} the lifted family is Zoll, so that $\tilde P^+$ and $\tilde P$ are closed manifolds. The parameter space $\tilde P^+$ and the bundle $\nu^+\colon E\to\tilde P^+$ are moreover orientable by Lemma~\ref{lem:lifted_orientations}, and we have a covering map $E \rightarrow \ST \tilde P^+$, which is a contactomorphism for $n+1 \geq 3$, by Proposition~\ref{prp:contacto_cover}. For $n+1=2$, we have $\tilde M\cong\tilde P^+\cong\Sph^2$ by Propositions~\ref{prp:even_rat_hom_sphere} and~\ref{prp:berger_even_dim_orient}, so that the covering map $\ST\tilde M\to\ST\tilde P^+$ is a covering map $\RP^3\to\RP^3$ and hence has degree one. Thus it is a contactomorphism in this case as well.

We denote the deck transformation of the covering $\tilde M \rightarrow M$ by $a$, and its lift to $E=\ST \tilde M$ by $\tilde a$. The latter acts freely, commutes with $\tau$ and preserves the coorientation of the contact structure.
Let $a^+$ and $\tau^+$ be the involutions of $\tilde P^+$ induced by $\tilde a$ and $\tau$. By Lemma~\ref{lem:contact_equiv}, $\tilde a$, $\tau$ and $\tilde a \tau$ correspond to $T a^+$, $-T\tau^+$ and $-T( a^+  \tau^+)$. The parameter space of the Besse family on $M$ is $P=\tilde P^+/K$, where $K$ is the group generated by $ a^+$ and $ \tau^+$. It is a manifold if $K$ acts freely.

We already know that $\tau^+$ is either the identity or that it acts freely. By Lemma~\ref{lem:fix} (a) the same is true for $a^+\tau^+$. Moreover, $a^+ \neq \mathrm{id}$ since $Ta^+$ acts freely. Together with Lemma~\ref{lem:fix} (b), we obtain the following consequence.

\begin{cor}\label{cor:fixed_a_plus} A Besse family on a closed non-orientable manifold $M$ is Zoll if and only if $a^+$ has no fixed point on $\tilde P^+$. If $a^+$ fixes a point $\sigma\in \tilde P^+$, then $T_\sigma a^+ = - \mathrm{id}$.
\end{cor}

The second claim of the next proposition will be used later in Section~\ref{sec:gromoll_grove}.

\begin{prp} \label{Berger_even_non_orientable}
Let $M$ be a closed non-orientable manifold of even dimension $n+1 \geq 2$. Then every Besse family on $M$ is Zoll, i.e. all its hypersurfaces are regular. Moreover, $\tau^+=a^+$.
\end{prp}
\begin{proof} The lift $\tilde a$ preserves the orientation of $E$ and it reverses the orientation of the fibers of $\nu^+$, which is defined by the coorientation of the contact structure, preserved by $\tilde a$, and the orientation of $\tilde M$, reversed by $a$. Therefore, $a^+$ reverses the orientation of $\tilde P^+$. But, by Corollary ~\ref{cor:fixed_a_plus}, at a fixed point $T_\sigma a^+=-\mathrm{id}$ preserves the orientation as $n+1$ is even. Hence $a^+$ has no fixed point and the Besse family on $M$ is Zoll by Corollary ~\ref{cor:fixed_a_plus}.

Since both $\tau^+$ and $a^+$ reverse the orientation of $\tilde P^+$, their product $\tau^+a^+$ preserves the orientation of $\tilde P^+$. Since $\tilde P^+$ is a rational homology sphere, the Lefschetz number of $\tau^+a^+$ satisfies $L(\tau^+a^+)=2$, and hence $\tau^+a^+$ has a fixed point on $\tilde P^+$ by the Lefschetz fixed point theorem \cite{Le26}; see also \cite[Chapter~3, \S 4]{GP74}, for instance. But then $\tau^+a^+$ is the identity by Lemma~\ref{lem:fix} (a) as already observed above, and so $\tau^+=a^+$ as $\tau^+$ and $a^+$ are involutions.
\end{proof} 

Now we deal with the odd-dimensional case.

\begin{prp} \label{Berger_odd_non_orientable}
Let $M$ be a closed non-orientable manifold of odd dimension $n+1 \geq 3$. Then every Besse family on $M$ is Zoll, i.e. all its hypersurfaces are regular.
\end{prp}
\begin{proof} Suppose again that $a^+$ fixes a point $\sigma\in \tilde P^+$. By Corollary~\ref{cor:fixed_a_plus}  all fixed points are isolated and at each fixed point the differential of $a^+$ is $-\mathrm{id}$. In particular, all fixed points are nondegenerate and have the same fixed point index $\mathrm{sign}\det(\mathrm{id}-T_\sigma a^+)=+1$, so by Lefschetz fixed point formula the Lefschetz number of $a^+$ satisfies $L(a^+)=|\mathrm{Fix}(a^+)|\geq 1$ \cite{Le26}. Moreover, $a^+$ reverses the orientation of $\tilde P^+$, because $\det(- \mathrm{id})=-1$ as the dimension is odd.

For an orientation-reversing involution $f$ of a closed oriented manifold $N$ of dimension $2m+1$ Poincar\'e duality gives $\mathrm{tr}(f^*|H^{2m+1-j}(N;\Q))=-\mathrm{tr}(f^*|H^{j}(N;\Q))$, and therefore
\[
        L(f) = \sum_{j=0}^{2m+1} (-1)^j \,\mathrm{tr}(f^* \mid H^j(N;\Q)) = 2\sum_{j=0}^{m} (-1)^j\, \mathrm{tr}(f^* \mid H^j(N;\Q))
\]
This applies to $a^+$ acting on $\tilde P^+$ and to $a$ acting on $\tilde M$. The $n$-sphere bundles $\mu^+\colon E\to\tilde M$ and $\nu^+\colon E\to\tilde P^+$ are orientable, since $\tilde M$ is orientable and by Lemma~\ref{lem:lifted_orientations}. Their Gysin sequences therefore give isomorphisms $H^j(\tilde M;\Q) \cong H^j(E;\Q)\cong H^j(\tilde P^+;\Q)$ for $j\leq n-1=2m-1$, which intertwine $a^*$, $\tilde a^*$ and $(a^+)^*$. Since $m\leq 2m-1$, the traces in the formula applied to $f=a$ and $f=a^+$ agree, so $L(a^+)=L(a)$. But $L(a)=0$, since $a$ acts freely, a contradiction. Hence, $a^+$ acts freely and so the Besse family on $M$ is Zoll by Corollary~\ref{cor:fixed_a_plus}.
\end{proof} 

\begin{proof}[Proof of Theorem~\ref{thm:Berger}]
The claim is shown in Propositions~\ref{prp:berger_even_dim_orient} and~\ref{Berger_odd_orientable} in the orientable case and in Propositions~\ref{Berger_even_non_orientable} and~\ref{Berger_odd_non_orientable} in the non-orientable case.
\end{proof}

\begin{rem}
For spheres $\Sph^{n+1}$ of even dimension at least six, Theorem~\ref{thm:Berger} also follows from the fact that $P^+$ is $(n-1)$-connected, together with \cite[Theorem~B]{LR23}. This was the argument in the first version of this paper.
\end{rem}

\section{Gromoll--Grove-type theorem for Besse families}\label{sec:gromoll_grove}

In this section, we prove Theorem~\ref{thm:Gromoll--Grove}, which states that the hypersurfaces of a Zoll family, in other words of an unoriented transitive family of closed hypersurfaces, on a manifold $M$ of even dimension $n+1 \geq 2$, must be embedded. By Proposition~\ref{prp:besse_impliese_cmpact} the manifold $M$ is closed, and hence so are $\ST M$ and $P^+$. We use this throughout this section without further mention.

We first assume that $M$ is orientable. In this case $M$ is a rational homology sphere without finite covers by Proposition~\ref{prp:even_rat_hom_sphere} and Corollary~\ref{cor:no_finite_covers}. Moreover, since $M$ is orientable and the dimension is even, the hypersurfaces in the family are $n$-spheres and $\nu$ is coorientable by Lemma~\ref{lem:lifted_orientations}, since all points of $P$ are regular. As in Section~\ref{sec:Berger_conjec} we consider the lifted sphere bundle
\[
  \nu^+ : \ST M \longrightarrow P^+ .
\]
By Proposition~\ref{prp:even_rat_hom_sphere} the parameter space $P^+$, which is a manifold, is a rational homology sphere as well and by Corollary~\ref{cor:no_finite_covers} it has no finite covers.

The proof of Theorem~\ref{thm:Gromoll--Grove} will be based on a contradiction to the nonvanishing of an Euler class. Since now $P^+$ is an orientable closed manifold, we deduce the following from Corollary~\ref{cor:euler_class}.

\begin{cor}\label{cor:euler}
$e(\nu^+)[P^+]=\pm2$. In particular, $e(\nu^+)\neq0$ in $H^{n+1}(P^+;\Z)$ and $H^{n+1}(P^+;\Q)$.
\end{cor}

\begin{proof}
Since $P^+$ is a manifold, the projection $BP^+\to P^+$ is a homotopy equivalence, and $H^{n+1}(P^+;\Z)\cong\Z$ since $P^+$ is closed, connected and orientable by Lemma~\ref{lem:lifted_orientations}. Hence the claim follows from $H^{n+1}(P^+;\Z)/\langle e(\nu^+)\rangle\cong\Z_2$; see Corollary~\ref{cor:euler_class}.
\end{proof}

Corollary~\ref{cor:euler} implies that $\nu^+$ admits no continuous section,
since a section would force $e(\nu^+) = 0$; cf. the proof of Proposition~\ref{prp:Gromoll_Grove_even_dimensional}. For the equatorial family, $P^+$ is the sphere of
poles and $\nu^+$ its unit tangent bundle, so this is the hairy ball
theorem. The statement can be applied in the proof of Theorem~\ref{thm:Gromoll--Grove} in the case $n=1$ (see Remark~\ref{rem:dim1_hairy_ball}) and underlies the idea of the proof in the general case: a double point would produce such a section, at least after pulling $\nu^+$ back to the \emph{space of ordered double points} that we define as
\[
  I := \big\{(z,w) \in \ST M \times \ST M \ :\ z \neq w,\
    \nu^+(z) = \nu^+(w),\ \mu^+(z) = \mu^+(w)\big\} \, .
\]
We next show that $I$ is a manifold and that 
$\pi_I \colon I \to P^+$, $(z,w) \mapsto \nu^+(z)$, defines a fiber bundle. Since its fibers have dimension less than $n$ and $P^+$ is a rational homology sphere without finite covers (see Proposition~\ref{prp:even_rat_hom_sphere} and Corollary~\ref{cor:no_finite_covers}), we will then see that pulling back along $\pi_I$ does not kill the top-dimensional class $e(\nu^+)$; see the proof of Theorem~\ref{thm:Gromoll--Grove}.

\begin{prp}\label{prop:I}
$I$ is a closed smooth $2n$-manifold and $\pi_I$ is a fiber bundle whose
fiber is a closed $(n-1)$-manifold, empty precisely when every hypersurface is
embedded. If $X$ is a connected component of $I$, then
$\pi_I|_X \colon X \to P^+$ is a fiber bundle with connected fiber.
\end{prp}

\begin{proof}
We first note that two branches of a hypersurface have distinct tangent hyperplanes at a double point. Indeed, if the two branches through a double point $(z,w)$ had the same tangent hyperplane, then $w=\tau(z)$. However, since $\nu$ is coorientable, the preimage in $\ST M$ of each fiber of $\nu$ consists of two fibers of $\nu^+$ interchanged by $\tau$, so $z$ and $\tau(z)$ never lie in a common fiber of $\nu^+$. This is where the family is used as an \emph{unoriented} transitive family.

We trivialize $\nu^+$ over a neighborhood $U$ of a given $\sigma_0 \in P^+$, writing a
point of $\ST M$ above $U$ as a pair $(\sigma, s)$ with $\sigma$ the fiber and
$s$ a point of it, and set $g(\sigma,s) := \mu^+(\sigma,s)$, so that
$g(\sigma,\cdot)$ is the immersion of $\Sigma_\sigma$. Writing $(z,w)=((\sigma,s),(\sigma,t))$, the triples $(\sigma,s,t)$ are local coordinates of the fiber product $\ST M\times_{P^+}\ST M=\{(z,w)\mid\nu^+(z)=\nu^+(w)\}$, which is a closed $(3n+1)$-dimensional submanifold of $\ST M\times\ST M$ since $\nu^+$ is a submersion. Near a double point $(\sigma_0,s_0,t_0)$, both $g(\sigma,s)$ and $g(\sigma,t)$ lie in a single chart $\kappa$ of $M$, and since $s_0\neq t_0$, the set $I$ coincides near $(\sigma_0,s_0,t_0)$ with the zero set of
\[
G(\sigma,s,t):=\kappa\, g(\sigma,s)-\kappa\, g(\sigma,t)\in\mathbb{R}^{n+1}.
\]

The partial differential of $G$ in $(s,t)$ has image $T\kappa(\Pi_1 + \Pi_2)$, where $\Pi_1,\Pi_2$ are the two branch tangent
hyperplanes. These are distinct by the previous paragraph and therefore
span $T_xM$, $x=g(\sigma_0,s_0)$. Hence $TG$ is surjective and $I$ is a submanifold of dimension
$(n+1)+n+n-(n+1) = 2n$. Moreover, $\pi_I$ is a submersion. Indeed, given
$\zeta \in T_{\sigma_0}P^+$, surjectivity in the fiber directions provides
$(v,u)$ with $TG(\zeta,v,u) = 0$, i.e.\ we have $(\zeta,v,u) \in T I$ and $T\pi_I (\zeta,v,u) = \zeta$ as needed. The fibers of $\pi_I$ have dimension $2n - (n+1)=n-1$.

We claim that $I$ is compact and that $\pi_I$ is proper. By compactness of $\ST M\times \ST M$ it suffices to show that $I$ is
closed. The conditions $\mu^+(z)=\mu^+(w)$ and $\nu^+(z)=\nu^+(w)$ are
closed, so it remains to exclude limits on the diagonal. If $(z_k,w_k)\in I$ converges to $(z,z)$, then for each $k$ the points $z_k\neq w_k$ have the same image under $(\mu^+,\nu^+)$, which is locally injective by Subsection~\ref{sub:double_Legendrian}, applied to the double Legendrian fibration $(\mu^+,\nu^+)$, a contradiction.

If $I\neq\emptyset$, then $\pi_I$ is surjective, since its image is open and closed in the connected space $P^+$. By Ehresmann's theorem \cite{Eh51}, $\pi_I$ is therefore a fiber bundle.

The fiber of $\pi_I$ over $\sigma \in P^+$ is identified with the set of ordered pairs
$(s,t)$ of distinct points of the fiber with $g(\sigma,s)=g(\sigma,t)$, that
is, with the ordered double points of $\Sigma_\sigma$. It is therefore empty
precisely when $\Sigma_\sigma$ is injectively immersed, which for a closed
hypersurface means embedded. Since $P^+$ is connected, all fibers are
diffeomorphic, so either every hypersurface is embedded or none is; cf. Lemma~\ref{lem:simple_charact}.

Finally, let $X$ be a connected component of $I$. It is open and closed, so $\pi_I|_X$ is a proper submersion, and its image is open and closed in the connected space $P^+$, hence all of it. Again by Ehresmann’s theorem, $\pi_I|_X$ is a fiber bundle. By compactness of its fibers $F$, the fibers have finitely many components. By connectedness of $X$, the fundamental group of $P^+$ acts transitively on the connected components of $F$. Since $\pi_1(P^+)$ has no proper finite index subgroup, this action, which factors through a finite symmetric group, must be trivial. Hence, $F$ is connected by transitivity of the action.
\end{proof}

Note that triple and higher multiplicity points cause no difficulty: each pair of branches meeting at such a point is separately transverse and contributes
its own point to $I$.

\begin{prp}\label{prp:Gromoll_Grove_even_dimensional} A Zoll family on an orientable manifold of even dimension $n+1 \geq 2$ is two-sided simple Zoll, i.e. all its hypersurfaces are two-sided and embedded.
\end{prp}
\begin{proof}
Suppose some hypersurface has a double point, so that $I \neq \emptyset$, and fix a
connected component $X$ of $I$. By Proposition~\ref{prop:I} the restriction
$\pi \colon X \to P^+$ of $\pi_I$ is a fiber bundle with connected fiber
$F$, a closed $(n-1)$-manifold.

The map $\mathrm{pr}_1 \colon X \to \ST M$, $(z,w)\mapsto z$, satisfies
$\nu^+\circ \mathrm{pr}_1 = \pi$, so $(z,w) \mapsto ((z,w),z)$ is a section $s$ of the pullback
$\varpi \colon \pi^*\ST M \to X$ of $\nu^+$, which is an orientable $\Sph^n$-bundle since $\nu^+$ is orientable.
From $\varpi\circ s = \mathrm{id}$ we get $s^*\circ \varpi^* = \mathrm{id}$, so $\varpi^*$
is injective in cohomology, and exactness of the Gysin sequence
\[H^0(X;\Q) \xrightarrow{\cup\, e(\varpi)} H^{n+1}(X;\Q) \xrightarrow{\varpi^*}
H^{n+1}(\pi^*\ST M;\Q)\] gives $e(\varpi) = 0$. By naturality of the Euler class,
\begin{equation}\label{eq:zero}
   \pi^*e(\nu^+) \;=\; e(\varpi) \;=\; 0 \qquad\text{in } H^{n+1}(X;\Q).
\end{equation}
Next we show that $\pi^*$ is injective in degree $n+1$, which contradicts
\eqref{eq:zero} because $e(\nu^+) \neq 0$ by Corollary~\ref{cor:euler}, and so
proves $I=\emptyset$, i.e. the embeddedness statement; see Proposition~\ref{prop:I}.

Consider the Serre spectral sequence of $F \to X \xrightarrow{\pi} P^+$.
The fiber is connected and the monodromy action on $H^q(F;\Q)$ defines a homomorphism from the finitely generated group $\pi_1(P^+)$ into $\mathrm{Aut}(H^q(F;\Q))$. Its image is a finitely generated linear group and hence residually finite \cite{Mal40}, so it is either trivial or has a proper subgroup of finite index. In the latter case, the preimage of this subgroup would be a proper subgroup of finite index in $\pi_1(P^+)$.
Since $\pi_1(P^+)$ has no proper subgroups of finite index, the monodromy is trivial and so we have
$E_2^{p,q} = H^p(P^+;H^q(F;\Q)) \cong  H^p(P^+;\Q) \otimes H^q(F;\Q)$. By Proposition~\ref{prp:even_rat_hom_sphere}
this vanishes unless $p \in \{0,n+1\}$, and it vanishes for $q \geq n$
because $\dim F = n-1$. Since no nontrivial differential can have target $E^{n+1,0}_r$, we have
$E_\infty^{n+1,0} = E_2^{n+1,0} = H^{n+1}(P^+;\Q)$. Moreover, $E_\infty^{n+1,0}$ is the only possibly nonzero term on its
antidiagonal: the other candidate is
$E_\infty^{0,n+1} \subseteq H^{n+1}(F;\Q) = 0$. So the filtration of
$H^{n+1}(X;\Q)$ has a single nonzero quotient and
$H^{n+1}(X;\Q) \cong E_\infty^{n+1,0}$. Since $\pi^*$ is the composition
$H^{n+1}(P^+;\Q) = E_2^{n+1,0} \twoheadrightarrow E_\infty^{n+1,0}
\hookrightarrow H^{n+1}(X;\Q)$ (see \cite[Theorem~5.9]{McC01}), and both maps are isomorphisms as seen above, $\pi^*$ is an isomorphism in degree $n+1$. 

Consequently, $I = \emptyset$ and so every hypersurface is embedded.

The hypersurfaces are two-sided since $\nu$ is coorientable, as noted at the beginning of this section; see Subsection~\ref{sub:double_Legendrian}.
\end{proof}

\begin{rem}\label{rem:dim1_hairy_ball}
For $n=1$ the fiber $F$ is a point, so $\pi\colon X\to P^+$ is a diffeomorphism, and $\mathrm{pr}_1\circ\pi^{-1}$ 
is a section of $\nu^+$, which cannot exist by Corollary~\ref{cor:euler}. In fact, $P^+\cong\Sph^2$ by Proposition~\ref{prp:even_rat_hom_sphere}, and $\nu^+$ is an oriented circle bundle over $\Sph^2$ with Euler number $\pm2$, i.e.\ isomorphic to the unit tangent bundle of $\Sph^2$ up to orientation, so this is again the hairy ball theorem.
\end{rem}

\begin{prp}\label{prp:Gromoll_Grove_non-orien} A Zoll family on a non-orientable manifold $M$ of even dimension $n+1 \geq 2$ is one-sided simple Zoll, i.e. all its hypersurfaces are one-sided and embedded.
\end{prp}
\begin{proof} 
Using the notations of Subsection~\ref{sub:Berger_non_orientable} we have $a^+=\tau^+$ by Proposition~\ref{Berger_even_non_orientable}. Hence $\tilde a$ maps each fiber $\Sigma^+$ of $\nu^+$ to $\tau(\Sigma^+)$, so that the deck transformation $a$ leaves every hypersurface of the lifted Zoll family on $\tilde M$ invariant and reverses its coorientation. Since $a$ acts transitively on the components of the preimage of each hypersurface of the Zoll family on $M$, every such hypersurface is the quotient of a hypersurface of the lifted family by $a$.
As $a$ acts freely on $\tilde M$ and the lifted hypersurfaces are embedded by Proposition~\ref{prp:Gromoll_Grove_even_dimensional}, the hypersurfaces of the Zoll family on $M$ are embedded and one-sided. In particular, $\nu$ is not coorientable, and so by Subsection~\ref{sub:coor_and_lifts}, the hypersurfaces of the Zoll family on $M$ are diffeomorphic to $\RP^n$.
\end{proof}

\begin{proof}[Proof of Theorem~\ref{thm:Gromoll--Grove}]
That $M$ is closed was shown in Proposition~\ref{prp:besse_impliese_cmpact}. The remaining claim in the orientable case is shown in Proposition~\ref{prp:Gromoll_Grove_even_dimensional} and in the non-orientable case it is shown in Proposition~\ref{prp:Gromoll_Grove_non-orien}.
\end{proof}

\section{A Wadsley-type theorem for minimal hypersurfaces}\label{sec:Wadsley}

In this section we prove Theorem~\ref{thm:Wadsley_intro}, an analogue of Wadsley's theorem for minimal hypersurfaces. We first recall the relevant background on flows and foliations with all leaves compact.

By Wadsley's theorem \cite{Wa75}, a smooth $\R$-action on a smooth manifold $M$ all of whose orbits are circles is orbit equivalent to a smooth circle action if and only if the orbits are totally geodesic with respect to some Riemannian metric on $M$. Its proof shows moreover that the flow of the unit vector field tangent to the orbits  with respect to such a Riemannian metric is periodic if $M$ is connected; cf.~\cite[Theorem~7.12]{Be78}. In particular, a Reeb flow, and hence the geodesic flow of a Riemannian manifold, all of whose orbits are closed is periodic. For geodesic flows on Riemannian orbifolds a Wadsley-type theorem was obtained in \cite{ALR21}.

Without a geometric assumption this conclusion fails. Sullivan constructed a flow on a compact $5$-manifold all of whose orbits are closed but whose periods are unbounded \cite{Su76}, and Epstein and Vogt constructed such a flow on a compact $4$-manifold \cite{EV78}. On compact $3$-manifolds, on the other hand, the periods of such flows are bounded \cite[Theorem~4.2]{Ep72}. More generally, one can ask whether a foliation of a compact manifold with all leaves compact has bounded leaf volume. This is equivalent to every leaf having finite holonomy \cite{Ep76}, that is, finite group of germs of return maps on a transverse disc along loops in the leaf; see \cite{CLN85,MM03} for background, for instance. This holds in codimension $2$, as shown independently by Edwards, Millett and Sullivan \cite{EMS77} and by Vogt \cite{Vo76}. In codimension at least $3$ it fails: products of the example of Epstein and Vogt with closed manifolds, foliated by points or by a single leaf, give compact foliations of compact manifolds with unbounded leaf volume in every codimension at least $3$ and every leaf dimension. With a geometric assumption, however, it holds in every codimension: by a theorem of Rummler, a foliation with all leaves compact of a compact manifold has bounded leaf volume if and only if its leaves are minimal submanifolds for some Riemannian metric \cite[Théorème~3]{Ru79}. Rummler points out that this generalizes Wadsley's theorem, except that Wadsley does not need compactness, and asks whether the implication from minimal leaves to locally bounded leaf volume holds on non-compact manifolds \cite[p.~236]{Ru79}.

The definition of a Besse family asks for an orbifold fiber bundle. For a foliation with all leaves compact, this amounts to every leaf having finite holonomy. Indeed, by the local Reeb stability theorem \cite[Theorem~2.9]{MM03}, every leaf then has a saturated neighborhood that is a finite quotient of a product, compatibly with the foliation. Hence the leaf space is an orbifold whose local groups are the holonomy groups \cite[Theorem~2.15]{MM03}, and the projection onto it is an orbifold fiber bundle; see also \cite[Theorem~4.3]{Ep76}. 

For an $(n+1)$-dimensional manifold $M$, $\PT M$ has dimension $2n+1$, and the Legendrian lifts of hypersurfaces of $M$ are $n$-dimensional, so a foliation of $\PT M$ by such lifts has codimension $n+1$. For $n=1$, where the hypersurfaces are curves, every foliation of $\PT M$ by their Legendrian lifts is an orbifold fiber bundle if $M$ is closed, by Epstein's theorem \cite{Ep72} (see Corollary~\ref{cor:Epstein_refinement}). For $n\geq2$ the codimension is at least $3$, where foliations with all leaves compact can have unbounded leaf volume, equivalently leaves with infinite holonomy. Theorem~\ref{thm:Wadsley} answers Rummler's question affirmatively for foliations of $\PT M$ by Legendrian lifts of minimal hypersurfaces: minimality of the hypersurfaces in $M$ suffices, in every dimension and without assuming $M$ to be compact. Throughout this section, $\mu\colon\PT M\to M$ denotes the projection and $\eta$ the natural contact structure on $\PT M$, as in Subsection~\ref{sub:Lege_fib_bun}. 

\begin{thm}\label{thm:Wadsley}
Let $(M,g)$ be a connected Riemannian manifold of dimension $n+1$ with $n\geq2$. Let $\mathcal F$ be a smooth foliation of $\PT M$ whose leaves $L$ are compact and are the Legendrian lifts of the immersions $\mu|_L\colon L\to M$. If every immersion $\mu|_L$ is minimal, then every leaf of $\mathcal F$ has finite holonomy. In particular, the leaf space $P=\PT M/\mathcal F$ is an orbifold, and the projection $\PT M\to P$ is an orbifold fiber bundle.
\end{thm}

\begin{proof}[Proof of Theorem~\ref{thm:Wadsley_intro}]
Since the Legendrian lifts are embeddings, each leaf $L$ of the foliation is the Legendrian lift of the immersion $\mu|_L$, which is a reparametrization of a hypersurface in the collection. Hence Theorem~\ref{thm:Wadsley} applies, and its conclusion is the definition of a Besse family.
\end{proof}

For closed $M$, Theorem~\ref{thm:Wadsley} also follows from \cite[Théorème~1]{Ru79} together with Lemma~\ref{lem:epstein} and its proof. Indeed, the leaves of $\mathcal F$ are minimal for the metric $\hat g$ from the proof of Lemma~\ref{lem:epstein}: since $\hat g\geq\mu^*g$ with equality on $T\mathcal F$, the $\hat g$-volume of a variation of a leaf $L$ is at least the $g$-area of its projection, with equality at $L$, so the first variation of the volume vanishes at $L$. Unlike Rummler's proof, ours is local and does not require $M$ to be compact.

Steps~1--3 of the proof of Theorem~\ref{thm:Wadsley} in Subsection~\ref{sub:proof_Wadsley} parallel the first half of Wadsley's argument \cite[\S4]{Wa75}. Step~4 replaces the second half by an argument specific to Legendrian foliations. The key object is the set $C_1$ of points near which the leaf area is unbounded (see Subsection~\ref{sub:compact_foliations}). It consists of the leaves with infinite holonomy and is nowhere dense (Lemmas~\ref{lem:epstein} and~\ref{lem:dense}). On each connected component $W$ of the complement of $C_1$, the holonomy is finite (Step~1), and the area of the leaves with trivial holonomy is constant by the first variation formula. Consequently the area is bounded on $W$ (Step~2). This step, the only one using minimality, parallels \cite[Corollary~4.4 and Proposition~4.5]{Wa75}. A counting argument going back to Epstein \cite[\S7]{Ep76} then shows that the linear holonomy of every leaf in the closure of $W$ is finite (Step~3). This is the analogue of \cite[Corollary~4.6]{Wa75}.

For the analogue of $W$, Wadsley then shows that its boundary orbits have locally bounded period, using a metric for which the flow is isometric on its closure \cite[Lemmas~4.8 and~4.9]{Wa75}. In our setting this part is replaced by a rigidity property of foliations by Legendrian lifts: the leaves carry an $\RP^n$-structure whose holonomy is the dual of the linear holonomy (Proposition~\ref{prp:proj}). For $n\geq2$, finite linear holonomy therefore forces finite fundamental group (Proposition~\ref{prp:dich}), hence finite holonomy. So the closure of $W$ contains no point of $C_1$, hence $W$ is closed, and $C_1=\emptyset$ (Step~4). For $n=1$ this argument fails, since $\pi_1(\RP^1)\cong\Z$ is infinite. There the leaves are circles, and the passage from finite linear holonomy to finite holonomy is what the second half of Wadsley's argument achieves.

In summary, minimality yields the area bound, and the Legendrian condition the passage from finite linear holonomy to finite holonomy.

\subsection{Compact foliations and leaf area}\label{sub:compact_foliations}

From now on, $(M,g)$ and $\mathcal F$ are as in Theorem~\ref{thm:Wadsley}. For a leaf $L$ we put
\[
\operatorname{area}(L)=\operatorname{area}_g(\mu|_L),
\]
the volume of $L$ for the leafwise metric $\mu^*g|_{T\mathcal F}$. Since the leaves are transverse to the fibers of $\mu$, this is a  Riemannian metric on $T\mathcal F$. We write $L_x$ for the leaf through $x$ and call a subset of $\PT M$ saturated if it is a union of leaves.

For a leaf $L$ and $p\in L$ we write $\Gamma=\pi_1(L,p)$, $h_\gamma$ for the germ at $p$ of the holonomy of $\gamma\in\Gamma$ on a transverse disc $Y$ through $p$, and $\rho(\gamma)=T_ph_\gamma$ for the linear holonomy, acting on the normal space $N_p=T_p\PT M/T_pL\cong T_pY$. The images of $\gamma\mapsto h_\gamma$ and $\gamma\mapsto\rho(\gamma)$ are the holonomy group and the linear holonomy group of $L$, respectively; see e.g.\ \cite[Chapter~2]{MM03}. Finally, let
\[
C_1=\{x\in\PT M\mid \operatorname{area}\text{ is unbounded on the leaves through every neighborhood of }x\}.
\]

\begin{lem}\label{lem:epstein}
Let $L$ be a leaf of $\mathcal F$. Then $\operatorname{area}$ is bounded on the leaves through some neighborhood of $L$ if and only if the holonomy group of $L$ is finite. Consequently, $C_1$ is the union of the leaves with infinite holonomy. In particular, it is closed and saturated.
\end{lem}

\begin{proof}
Since the leaves are transverse to the fibers of $\mu$, we have $T\mathcal F\cap\ker T\mu=0$, so there is a complement $H\supset T\mathcal F$ of $\ker T\mu$ in $T\PT M$. Let $\mathrm{pr}\colon T\PT M\to\ker T\mu$ be the projection along $H$, choose a Riemannian metric $g_V$ on $\ker T\mu$, and put $\hat g=\mu^*g+\mathrm{pr}^*g_V$. This is a Riemannian metric on $\PT M$ with $\hat g\geq\mu^*g$ and $\hat g|_{T\mathcal F}=\mu^*g|_{T\mathcal F}$. Hence $\operatorname{area}(L)$ is the Riemannian volume of $L$ for $\hat g$. By the construction in \cite[3.5.2 and 3.5.3]{Ep76}, this is a volume function in the sense of \cite[3.1]{Ep76}. The first claim is therefore the equivalence of (4.2.1) and (4.2.3) in \cite[Theorem~4.2]{Ep76}. Since saturations of open sets are open, the area is bounded on the leaves through some neighborhood of a point $x$ if and only if it is bounded on the leaves through some neighborhood of $L_x$. Hence $C_1$ is the union of the leaves with infinite holonomy. It is closed by definition and saturated as a union of leaves.
\end{proof}

The next lemma is the analogue of Wadsley's observation that the points of continuity of the period function are dense \cite[p.~545]{Wa75}.

\begin{lem}\label{lem:dense}
The set $C_1$ is nowhere dense, i.e.\ $\PT M\setminus C_1$ is dense in $\PT M$.
\end{lem}

\begin{proof}
The area is lower semicontinuous at compact leaves \cite[\S9, p.~281]{Ep76}, so the sets $Q_k=\{x\mid\operatorname{area}(L_x)\leq k\}$, $k\in\N$, are closed. They cover $\PT M$ since leaves are compact. Let $U\subset\PT M$ be open and nonempty. By the Baire category theorem applied to $U=\bigcup_k(Q_k\cap U)$, some $Q_k\cap U$ contains a nonempty open set $U'$. The area is bounded by $k$ on the leaves through $U'$, so $U'\cap C_1=\emptyset$. Hence $\PT M\setminus C_1$ is dense, and since $C_1$ is closed, it is nowhere dense.
\end{proof}

\subsection{A pointwise counting lemma}

The next lemma is a pointwise version of two steps in Epstein's proof that bounded volume implies finite holonomy: the stabilizing chain of iterated holonomy images in the proof of \cite[Proposition~7.1]{Ep76}, and the passage from bounded orbits to a finite-index subgroup acting trivially in \cite[\S7.2]{Ep76}. Epstein applies these steps on a saturated neighborhood of $L$. In the proof of Theorem~\ref{thm:Wadsley}, however, a bound is available only for the leaves through an open subset of a transverse disc at $p$ that accumulates at $p$ without being a neighborhood of $p$. As we have not found such a pointwise statement in the literature, we include a proof.

Let $L$ be a leaf, $p\in L$ and $Y$ a transverse disc at $p$. Let $S$ be a finite generating set of $\Gamma=\pi_1(L,p)$, which exists since $L$ is compact, let $F$ be the free group on $S$, and let $\pi\colon F\to\Gamma$ be the canonical surjection. By a word we mean a finite sequence of letters $s^{\pm1}$, $s\in S$. Each word $w$ represents an element of $F$, and $|w|$ denotes its length. For each $s\in S$ fix a loop $c_s$ at $p$ representing $s$ and a chain of foliated charts covering it, and let $h_s\colon D_s\to Y$ be the holonomy transport along this chain, defined on a neighborhood $D_s$ of $p$ in $Y$. Let $h_{s^{-1}}=h_s^{-1}\colon h_s(D_s)\to D_s$ be the transport along the reversed chain. These maps send each point $x$ of their domain to a point of $L_x\cap Y$. For a word $w=s_1^{\varepsilon_1}\cdots s_m^{\varepsilon_m}$ put $h_w=h_{s_1^{\varepsilon_1}}\circ\dots\circ h_{s_m^{\varepsilon_m}}$ on its natural domain. Then $h_w(x)\in L_x\cap Y$ whenever $h_w(x)$ is defined, and $h_w$ represents the germ $h_{\pi(w)}$, so $T_ph_w=\rho(\pi(w))$. For $N\geq1$ let $F_N\subset F$ be the intersection of all subgroups of index at most $N$. Since $F$ is finitely generated, it has only finitely many such subgroups, as these correspond to transitive actions on finite sets. Hence $F_N$ has finite index in $F$, so it is finitely generated, and $\pi(F_N)$ has finite index in $\Gamma$.

\begin{lem}\label{lem:orbits}
Let $N\geq1$ and fix words $w_1,\dots,w_k$ that generate $F_N$. Then there is a neighborhood $V$ of $p$ in $Y$ with the following property: if $x\in V$ and $\#(L_x\cap Y)\leq N$, then $h_{w_i}(x)=x$ for $i=1,\dots,k$.
\end{lem}

\begin{proof}
Since there are finitely many words of length at most $N$, there is a neighborhood $V$ of $p$ in $Y$ on which $h_w$ is defined for all words $w$ of length at most $N$.

Let $x\in V$ with $\#(L_x\cap Y)\leq N$, and for $j\leq N$ put $\Omega_j=\{h_w(x)\mid |w|\leq j\}\subset L_x\cap Y$. These sets increase with $j$, start with $\Omega_0=\{x\}$, and have at most $N$ elements, so there is $j_0\leq N-1$ with $\Omega_{j_0}=\Omega_{j_0+1}=:\Omega$. If $y=h_w(x)\in\Omega$ with $|w|\leq j_0$, then $|s^{\pm1}w|\leq j_0+1\leq N$, so $h_{s^{\pm1}}(y)=h_{s^{\pm1}w}(x)$ is defined by the choice of $V$, and it lies in $\Omega_{j_0+1}=\Omega$. Hence $h_s$ and $h_{s^{-1}}$ are defined on $\Omega$ and map it into itself, and since $h_{s^{-1}}=h_s^{-1}$, their restrictions to $\Omega$ are mutually inverse bijections. So $s\mapsto h_s|_\Omega$ defines an action of $F$ on $\Omega$, in which a word $w$ acts by $h_w|_\Omega$. The orbit of $x$ has at most $N$ elements, so its stabilizer has index at most $N$ and contains $F_N$. Hence $h_{w_i}(x)=x$ for $i=1,\ldots,k$.
\end{proof}

\subsection{Leaves with finite holonomy}\label{sub:finite_holonomy}

We recall two standard facts about foliations with compact leaves. They hold for every smooth foliation whose leaves are compact. We apply them to the foliation $\mathcal F$ in Theorem~\ref{thm:Wadsley}. Lemma~\ref{lem:model} describes $\mathcal F$ near a leaf with finite holonomy, and Lemma~\ref{lem:orbifold_leafspace} describes the leaf space of a saturated region in which all leaves have finite holonomy.

\begin{lem}\label{lem:model}
Let $L$ be a leaf with finite holonomy group $H$, and let $\hat L\to L$ be the covering corresponding to the kernel of the holonomy homomorphism $\pi_1(L)\to H$. Then there are a transverse disc $T$ at a point $t_0$ of $L$, on which $H$ acts by diffeomorphisms fixing $t_0$, and a leaf-preserving diffeomorphism $\Psi$ from $\hat L\times_HT$ onto a saturated neighborhood of $L$. In particular, the leaf through $\Psi([\hat l,t])$ is $\hat L/\mathrm{Stab}_H(t)$ with holonomy group $\mathrm{Stab}_H(t)$, and as $t\to t_0$ the maps $\hat L\to\PT M$, $\hat l\mapsto\Psi([\hat l,t])$, converge in $C^1$ to the covering $\hat L\to L$.
\end{lem}

\begin{proof}
This is the local Reeb stability theorem \cite[Theorem~2.9]{MM03}. The remaining statements follow directly from the description of the saturated neighborhood as $\hat L\times_HT$.
\end{proof}

\begin{lem}\label{lem:orbifold_leafspace}
Let $U\subset\PT M$ be open and saturated, and suppose that every leaf in $U$ has finite holonomy. Then $U/\mathcal F$ is an orbifold, its local group at a leaf is the holonomy group of that leaf, its regular points are the leaves with trivial holonomy, and the projection $U\to U/\mathcal F$ is an orbifold fiber bundle.
\end{lem}

\begin{proof}
The first two statements are \cite[Theorem~2.15]{MM03}, applied to the restriction of $\mathcal F$ to $U$. The others follow from the description of saturated neighborhoods of leaves in Lemma~\ref{lem:model}, since $H$ acts effectively on $T$, the holonomy group being a group of germs.
\end{proof}

\subsection{The projective structure on Legendrian leaves}\label{sub:projective_structure}

The next lemma is a local version of the map $\Phi$ of Proposition~\ref{prp:contacto_cover}. For foliations, the global covering property of $\Phi$ is replaced by Propositions~\ref{prp:proj} and~\ref{prp:dich} (cf.\ also \cite[Theorem~6.2 and Corollary~6.3]{Pa90}).

\begin{lem}\label{lem:normalform}
Let $U\subset\PT M$ be open and let $\nu_U\colon U\to B$ be a submersion onto an open subset $B\subset\R^{n+1}$, with $\ker T\nu_U=T\mathcal F$. Then
\[
\Phi\colon U\to\PT B,\qquad z\mapsto T\nu_U(\eta_z),
\]
is a local contactomorphism. In particular, its restriction to a fiber $\nu_U^{-1}(y)$ is a local diffeomorphism onto an open subset of $\mathbb{P}(T^*_yB)\cong\RP^n$.
\end{lem}

\begin{proof}
Since the leaves are Legendrian, $\ker T\nu_U=T\mathcal F\subset\eta$. The computation in Subsection~\ref{sub:Lege_fib_bun} is pointwise and only uses this inclusion, so it shows that $\Phi$ is a local contactomorphism. The map $\Phi$ sends the fiber $\nu_U^{-1}(y)$ into the fiber $\mathbb{P}(T^*_yB)$. Both have dimension $n$, and the differential of $\Phi|_{\nu_U^{-1}(y)}$ is the restriction of the injective map $T\Phi$ to the tangent spaces of $\nu_U^{-1}(y)$. Hence $\Phi|_{\nu_U^{-1}(y)}$ is a local diffeomorphism.
\end{proof}

The following proposition constructs an $\RP^n$-structure, that is, a $(\mathrm{PGL}(n+1,\R),\RP^n)$-structure, on each leaf (cf.\ \cite[Chapter~3]{Th97} for $(G,X)$-structures, developing maps and holonomy).

\begin{prp}\label{prp:proj}
Every leaf $L$ carries an $\RP^n$-structure whose holonomy group is the image of the linear holonomy group of $L$ under the homomorphism $\mathrm{GL}(N_p)\to\mathrm{PGL}(N_p^*)$, $A\mapsto\mathbb P\bigl((A^{-1})^T\bigr)$, for $p\in L$.
\end{prp}
\begin{proof}
We cover $\PT M$ by foliated charts $U$ with projections $\nu_U\colon U\to B$ as in Lemma~\ref{lem:normalform}, chosen so small that $\Phi$ is injective on each fiber of $\nu_U$. By Lemma~\ref{lem:normalform}, the restrictions of the maps $\Phi$ to the plaques of $L$, i.e.\ the fibers of the maps $\nu_U$ contained in $L$, are diffeomorphisms onto open subsets of $\mathbb P(T^*_yB)\cong\RP^n$, and they cover $L$. Near each point of an overlap $U\cap U'$, the submersions $\nu_U$ and $\nu_{U'}$ are constant on the same plaques, so $\nu_{U'}=\psi\circ\nu_U$ for a diffeomorphism $\psi$ between open subsets of $B$ and $B'$. Hence $\Phi'(z)=T_y\psi\bigl(\Phi(z)\bigr)$ for $z\in\nu_U^{-1}(y)$.  A hyperplane $\zeta$ corresponds to the line of covectors vanishing on it, and the linear map $T^*_yB\to T^*_{\psi(y)}B'$, $\alpha\mapsto\alpha\circ(T_y\psi)^{-1}$, maps the annihilator of $\zeta$ to that of $T_y\psi(\zeta)$. Hence the transition maps are locally the projective transformations $\mathbb P\bigl(((T_y\psi)^{-1})^T\bigr)$. Since a projective transformation is determined by its restriction to an open set, they are projective on each connected component of an overlap of plaques. So these restrictions form an $\RP^n$-atlas.

For the holonomy, identify $N_p$ with $T_{y_0}B$ via $T_p\nu_U$ for a chart $U$ around $p$, with $y_0=\nu_U(p)$. Continuing charts along a loop $\gamma$ at $p$ composes the transition maps $\mathbb P\bigl(((T\psi_j)^{-1})^T\bigr)$. Since $A\mapsto(A^{-1})^T$ is a homomorphism and the composition of the $\psi_j$ represents the holonomy germ of $\gamma$, the holonomy of the $\RP^n$-structure along $\gamma$ is $\mathbb P\bigl((\rho(\gamma)^{-1})^T\bigr)$.
\end{proof}

The proof of the next proposition is a standard covering argument for geometric structures with compact model space.

\begin{prp}\label{prp:dich}
Let $n\geq2$ and let $L$ be a leaf whose linear holonomy group $\rho(\Gamma)$ is finite. Then $\pi_1(L)$ is finite, and $L$ has finite holonomy.
\end{prp}

\begin{proof}
By Proposition~\ref{prp:proj}, the holonomy group $G$ of the $\RP^n$-structure on $L$ is a homomorphic image of $\rho(\Gamma)$, so it is finite. Let $\mathrm{hol}\colon\Gamma\to G$ be the holonomy, $K\subset\Gamma$ its kernel and $L_K\to L$ the covering corresponding to $K$. It has $|G|$ sheets, so $L_K$ is compact, as $L$ is. A developing map $\tilde L\to\RP^n$ of the $\RP^n$-structure is $\mathrm{hol}$-equivariant, so it descends to a local diffeomorphism $L_K\to\RP^n$, which is a covering map since $L_K$ is compact. As $n\geq2$, we have $\pi_1(\RP^n)\cong\Z_2$, so $|K|\leq2$ and $|\pi_1(L)|\leq2|G|$. The second claim follows, as the holonomy group of $L$ is a quotient of $\Gamma$.
\end{proof}

\subsection{Proof of Theorem~\ref{thm:Wadsley}}\label{sub:proof_Wadsley}
\emph{Step 1.} By Lemmas~\ref{lem:epstein} and~\ref{lem:dense}, $C_1$ is closed and saturated, and its complement is nonempty. Let $W$ be a connected component of $\PT M\setminus C_1$. It is open, since the components of the open set $\PT M\setminus C_1$ are open, and saturated, since $C_1$ is saturated and the leaves are connected. Since $\operatorname{area}$ is locally bounded on $W$, every leaf in $W$ has finite holonomy by Lemma~\ref{lem:epstein}. We show in Steps~2--4 that $W$ is also closed. Since $\PT M$ is connected, this implies $W=\PT M$ and $C_1=\emptyset$.

\emph{Step 2.} Let $W_{\mathrm{reg}}\subset W$ be the union of the leaves with trivial holonomy. By Step~1 and Lemma~\ref{lem:orbifold_leafspace}, $W/\mathcal F$ is an orbifold with regular part $W_{\mathrm{reg}}/\mathcal F$. Since $W$ is connected, so is $W/\mathcal F$, and hence so is $W_{\mathrm{reg}}/\mathcal F$ by \cite[Lemma~5.5]{La18}. The projection $W_{\mathrm{reg}}\to W_{\mathrm{reg}}/\mathcal F$ is open with connected fibers, so $W_{\mathrm{reg}}$ is connected.

Near a leaf $L_0\subset W_{\mathrm{reg}}$, the foliation is a product $L_0\times Y$ with $Y$ a disc, and the restrictions of $\mu$ to the leaves form a smooth family of minimal immersions of the closed manifold $L_0$. By the first variation formula, the area is constant in this family. So the area is locally constant on $W_{\mathrm{reg}}$, hence equal to a constant $A_W$.

Since the regular part of an orbifold is dense and the projection $W\to W/\mathcal F$ is open, $W_{\mathrm{reg}}$ is dense in $W$. By Lemma~\ref{lem:model}, every leaf in $W_{\mathrm{reg}}$ near a leaf $L'\subset W$ is the image of an embedding of the covering $\hat L'$ of $L'$ from that lemma, and these embeddings converge in $C^1$ to this covering. Since the area is continuous under $C^1$ convergence and a covering of $L'$ has area at least $\operatorname{area}(L')$, every leaf in $W$ has area at most $A_W$.

\emph{Step 3.} Let $p\in\overline W$ and $L=L_p$. Choose a foliated chart $U\cong Z\times Y_0$ around $p$, with plaques $Z\times\{y\}$, and let $Y=\{a\}\times Y_0$ be the transverse disc through $p$ given by a slice of this chart. After shrinking the chart, the plaques have area at least some $\delta>0$ (compare \cite[2.3]{Ep76}). Since they are pairwise disjoint, $\operatorname{area}(L')\geq\delta\cdot\#(L'\cap Y)$ for every leaf $L'$. Since $W$ is saturated, $W\cap U$ is the union of the plaques through $W\cap Y$. As $W$ is open and $p\in\overline W$, the set $W\cap Y$ is open in $Y$ and has $p$ in its closure. Every leaf through $W\cap Y$ lies in $W$, so by Step~2 it has area at most $A_W$ and meets $Y$ in at most $N=\lfloor A_W/\delta\rfloor \geq 1$ points. With $S$, $F$, $F_N$, $\pi$ and $h_w$ as before Lemma~\ref{lem:orbits}, fix words $w_1,\dots,w_k$ generating $F_N$. By Lemma~\ref{lem:orbits}, and since the natural domain of each $h_{w_i}$ is open and contains $p$, there is a neighborhood $V$ of $p$ in $Y$, contained in all these domains, such that $h_{w_i}=\mathrm{id}$ on $W\cap V$ for all $i$. By continuity of the derivative,
\[
\rho(\pi(w_i))=T_ph_{w_i}=\lim_{x\to p,\ x\in W\cap V}T_xh_{w_i}=\mathrm{id}.
\]
Hence $\rho$ is trivial on the finite-index subgroup $\pi(F_N)$ of $\Gamma$, and $\rho(\Gamma)$ is finite.

\emph{Step 4.} By Proposition~\ref{prp:dich}, which uses $n\geq2$, the leaf $L$ has finite holonomy, so $p\notin C_1$ by Lemma~\ref{lem:epstein}. Since $\PT M\setminus C_1$ is open and $p\in\overline W$, the component of $\PT M\setminus C_1$ containing $p$ meets $W$, so $p\in W$. Hence $W$ is closed, and $C_1=\emptyset$ by Step~1. Thus every leaf has finite holonomy, and the remaining assertions follow from Lemma~\ref{lem:orbifold_leafspace} with $U=\PT M$. \qed
\newline

By Proposition~\ref{prp:dich} and its proof, every leaf then has finite fundamental group, and the leaves with trivial holonomy are diffeomorphic to $\Sph^n$ or $\RP^n$, since for them $G$ is trivial and $L$ itself covers $\RP^n$.

\section{Topological classification of Zoll families of embedded hypersurfaces}\label{sec:topo_class}

A surface with a Zoll Riemannian metric is diffeomorphic to $\mathbb{S}^2$ or $\mathbb{RP}^2$. In higher dimensions Zoll Riemannian metrics exist on all compact rank one symmetric spaces, i.e.\ spheres, real, complex and quaternionic projective spaces and the Cayley plane. In fact, by a theorem of Bott and Samelson \cite{Bo54, Sa63} (see also \cite[Theorem~7.23]{Be78}) every Zoll manifold has the integral cohomology ring of one of them. For Zoll families of hypersurfaces the situation is more rigid, as the following result shows. For $n\geq 2$, a stronger version is due to Martins \cite{Ma26}. We give a short proof of Theorem~\ref{thm:top_class} that covers all $n\geq 1$.

\begin{thm}\label{thm:top_class} For $n\geq 1$ let $\mathscr{Z}$ be an unoriented transitive family of closed embedded hypersurfaces of an $(n+1)$-manifold $M$. 

(i) If the hypersurfaces in $\mathscr{Z}$ are two-sided, then they are diffeomorphic to $\Sph^{n}$, the manifold $M$ is homeomorphic to $\Sph^{n+1}$, and the parameter space $P$ is homeomorphic to $\RP^{n+1}$. 

(ii) If the hypersurfaces in $\mathscr{Z}$ are one-sided, then they are diffeomorphic to $\RP^{n}$, and the manifold $M$ as well as the parameter space $P$ are homeomorphic to $\RP^{n+1}$.
\end{thm}

Note that by Lemma~\ref{lem:simple_charact}, it suffices to assume that one hypersurface of $\mathscr Z$ is embedded.

The theorem could equivalently be phrased in terms of a simple Zoll double Legendrian fibration for which at least one of the two bundles is not coorientable. Coorientability of the other bundle then determines whether we are in case~(i) or~(ii); see Subsection~\ref{sub:double_Legendrian}. 
If, instead, both bundles are coorientable, then an arbitrary base manifold can occur, though the fibers are still spheres; see Subsection~\ref{sub:coor_and_lifts} for an example, and \cite{Gu26} for examples by minimal hypersurfaces on simply connected compact rank one symmetric spaces.

For $n\geq 2$, \cite{Ma26} shows more strongly that $M$ is diffeomorphic to $\Sph^{n+1}$ resp.\ $\RP^{n+1}$, and in case~(ii) that $P$ is diffeomorphic to $\RP^{n+1}$. For $n=1$ and families consisting of the geodesics of a projective structure, Theorem~\ref{thm:top_class} is contained in \cite[Lemma~2.8, Propositions~2.9 and~2.11, and Corollaries~2.19 and~2.20]{LM02}; cf.\ \cite[Remark~2.9]{Ma26}. For $n+1\leq 3$, all homeomorphisms in Theorem~\ref{thm:top_class} can be replaced by diffeomorphisms, since homeomorphic manifolds of dimension at most three are diffeomorphic; see \cite{Th97}, for instance.

As explained in the introduction, Martins' proof and ours rest on the same core idea of isotoping a two-sided hypersurface of the family onto itself with reversed coorientation through other members of the family. Our proof is shorter and, apart from the argument of \cite[Proposition~3.5]{Ma26} (see below), independent. Moreover, it isolates the topological mechanism behind the argument; see Subsection~\ref{sub:reversible_submanifolds}. 

A first observation, that all hypersurfaces are diffeomorphic to $\Sph^n$ or to $\RP^n$, follows from the fact, well known in contact geometry, that $\Phi\colon\PT M\to\PT P$ is a local contactomorphism \cite{Ar89, Pa90}, together with the closedness of the hypersurfaces; see Section~\ref{sec:contact_geometry}. In \cite{Ma26} the local diffeomorphism property is obtained by a computation using an auxiliary Riemannian metric. In carrying out the core idea, we do not need to analyze nodal sets as in \cite{Ma26}. In contrast to \cite{Ma26}, our argument yields the classification only up to homeomorphism, since we identify both $M$ and $P$ via the generalized Poincaré conjecture. In \cite{Ma26} it is used only to identify $P$ in case~(i), while the diffeomorphisms are obtained via explicit isotopy and gluing constructions. To identify the relevant quotients of spheres by free involutions with $\RP^{n+1}$, we follow the argument of \cite[Proposition~3.5]{Ma26}.

Hypersurfaces that can be isotoped onto themselves with reversed coorientation, which we call \emph{isotopically reversible}, are the central notion of this section. We study them in greater generality than needed for Theorem~\ref{thm:top_class}. In particular, a manifold containing a hypersurface homeomorphic to a sphere that is isotopically reversible by an isotopy fixing a point is itself homeomorphic to a sphere; see Lemma~\ref{lem:sphere_fibers_implies_sphere_total_space}.

\subsection{Isotopically reversible submanifolds}\label{sub:reversible_submanifolds} 

Throughout this subsection, submanifolds are assumed to be
embedded and closed as subsets of the ambient manifold.

The following condition will come up naturally in our discussion of unoriented transitive families of closed embedded hypersurfaces.

\begin{dfn}A connected two-sided codimension one submanifold $\Sigma$ of a manifold $M$ is called \emph{isotopically reversible} if there exists a diffeotopy $\Psi:[0,1] \times M \rightarrow M$ of $M$ starting at the identity whose time-$1$ map $\Psi_1$ preserves $\Sigma$ and reverses the orientation of its normal bundle. We call such a $\Psi$ a \emph{reversing diffeotopy}. The submanifold is called \emph{pointed isotopically reversible} if there is a reversing diffeotopy $\Psi$ that in addition fixes a point on $\Sigma$.
\end{dfn}

For example, for any closed manifold $N$, the product $\Sph^1 \times N$ is pointed isotopically reversible in $\Sph^2 \times N$, where $\Sph^1$ sits inside $\Sph^2$ as an equator. In particular, $\Sph^1 \times \Sph^1$ is pointed isotopically reversible inside $\Sph^2 \times\Sph^1$. In these examples the hypersurface is not only isotopically reversible but also \emph{separating}, i.e.\ its complement has two connected components. This is a consequence of the following lemma.

\begin{lem}\label{lem:suf_cond_separating} Suppose a submanifold $\Sigma$ of a manifold $M$ is isotopically reversible. Then $\Sigma$ is separating in $M$.
\end{lem}
\begin{proof}
We choose an embedding $[-1,1]\times\Sigma \rightarrow M$ with closed image that restricts to the inclusion on $\{0\}\times\Sigma$, i.e.\ a closed tubular neighborhood of $\Sigma$, and define a continuous map from $M$ to $[-1,1]/(-1\sim 1)\cong\Sph^1$ as follows: on the tubular neighborhood we first project to $[-1,1]$ and then to $\Sph^1$. Everything outside the tubular neighborhood is collapsed to  the equivalence class of $1$ in $\Sph^1$. Pulling back the generator of $H^1(\Sph^1;\Z)$ via this map we obtain a class $\varphi \in H^1(M;\Z)$ which counts the signed intersection number of a loop in $M$ with $ \Sigma$. Suppose now that $ \Sigma$ is not separating. Then there is a loop $\gamma$ in $M$ that intersects $ \Sigma$ only once and transversally, so that $|\varphi([\gamma])|=1$. 
After reversing $\gamma$ if necessary, we can assume that $\varphi([\gamma])=1$. Let $\Psi_1$ be the time-$1$ map of a reversing diffeotopy. Since it flips the orientation of the normal bundle, we have $\varphi([\Psi_1\circ\gamma])=-1$. On the other hand, homotopy invariance implies $\varphi([\Psi_1\circ\gamma])=1$, a contradiction.
\end{proof}

Our standing assumption that $\Sigma$ is closed as a subset of $M$ rules out examples like $\Sigma=(-1,1)\times \{0\} \subset \R^2$, which are isotopically reversible but not separating. If the isotopically reversible hypersurface is in fact closed as a manifold, then also the ambient manifold has to be closed, as the following lemma shows.

\begin{lem}\label{lem:closed} Suppose an  isotopically reversible submanifold $\Sigma$ of a manifold $M$ is closed as a manifold. Then also $M$ is a closed manifold.
\end{lem}
\begin{proof}
By Lemma~\ref{lem:suf_cond_separating}, $\Sigma$ separates $M$ into two connected components $U_1$ and $U_2$. Let $\Psi$ be a reversing diffeotopy. Since $\Psi_1$ preserves $\Sigma$ and reverses the orientation of its normal bundle, it interchanges $U_1$ and $U_2$. We set
\[ D:=\{\,y\in M \mid \Psi_t(y)\in \Sigma \text{ for some } t\in[0,1]\,\}. \]
The map $\Sigma\times[0,1]\to M$, $(s,t)\mapsto \Psi_t^{-1}(s)$, is continuous, and hence its image $D$ is compact. Clearly, $\Sigma\subset D$. For $y\in U_1$ the path $t\mapsto\Psi_t(y)$ connects $y$ to $\Psi_1(y)\in U_2$ and hence meets $\Sigma$, so $U_1\subset D$. Likewise $U_2\subset D$. Hence $D=M$ is compact.
\end{proof}

Being isotopically reversible has further strong homological implications, as the following lemma shows.

\begin{lem}\label{lem:surjection_on_homology} Let $\Sigma$ be an isotopically reversible  submanifold of a manifold $M$ and let $\bar U_1$ and $ \bar U_2$ be the closures of the connected components of the complement of $\Sigma$ in $M$ (which exist by Lemma~\ref{lem:suf_cond_separating}). Then the inclusion-induced maps $H_*(\Sigma) \rightarrow H_*(\bar U_1)$ and $H_*(\Sigma) \rightarrow H_*(\bar U_2)$ are surjective.
\end{lem}
\begin{proof}
Let $\Psi$ be a reversing diffeotopy. As in the proof of Lemma~\ref{lem:closed}, $\psi=\Psi_1$ interchanges $\bar U_1$ and $\bar U_2$. The Mayer--Vietoris sequence reads
\[
\begin{tikzcd}[column sep=large]
\ldots \arrow[r] &  H_{k}(\Sigma) \arrow[r,"{((i_1)_*,\,(i_2)_*)}"] & H_k(\bar U_1) \oplus H_k(\bar U_2) \arrow[r,"{(j_1)_*-(j_2)_*}"] & H_k(M) \arrow[r] & \ldots 
\end{tikzcd}
\]
Since each cycle $z$ in $\bar U_1$ is homologous in $M$ to the cycle $\psi_*(z)$ in $\bar U_2$, the kernel of $(j_1)_*-(j_2)_*$ contains the graph of $\psi_*:H_k(\bar U_1)\rightarrow H_k(\bar U_2)$. By exactness of the sequence the map $(i_1)_*: H_k(\Sigma) \rightarrow H_k(\bar U_1)$ is therefore surjective. By symmetry also $(i_2)_*$ is surjective.
\end{proof}

For the sake of a simpler proof we state the following analogue of Lemma~\ref{lem:surjection_on_homology} on the level of fundamental groups only for a simply connected hypersurface $\Sigma$.

\begin{lem}\label{lem:surjection_on_pi1} Let $\Sigma$ be a simply connected, pointed isotopically reversible  submanifold of a manifold $M$ and let $\bar U_1$ and $ \bar U_2$ be the closures of the connected components of the complement of $\Sigma$ in $M$ (which exist by Lemma~\ref{lem:suf_cond_separating}). 
Then also $\bar U_1$, $\bar U_2$ and $M$ are simply connected.

\end{lem}
\begin{proof} We choose a basepoint $p$ in $\Sigma$ and a reversing diffeotopy $\Psi$ fixing $p$. By the Seifert--van Kampen theorem the fundamental group of $M$ is isomorphic to the free product
$\pi_1(\bar U_1) * \pi_1(\bar U_2)$. Since $\Psi$ fixes $p$, every loop in $\bar U_1$ based at $p$ is homotopic via $\Psi$ to a loop in $\bar U_2$ based at $p$, and vice versa. Hence the images of $\pi_1(\bar U_1)$ and $\pi_1(\bar U_2)$ in the free product coincide. On the other hand, in a free product the images of the two factors intersect trivially, so if they coincide, both are trivial, and hence so is $\pi_1(M)$.
\end{proof}

Combining the previous two lemmas with the generalized Poincar\'e conjecture, we can show the following statement.

\begin{lem}\label{lem:sphere_fibers_implies_sphere_total_space} Suppose a pointed isotopically reversible submanifold $\Sigma$ of an $(n+1)$-manifold $M$, $n\geq 1$, is homeomorphic to an $n$-sphere. Then $M$ is homeomorphic to an $(n+1)$-sphere.
\end{lem}
\begin{proof} By Lemma~\ref{lem:suf_cond_separating}, $\Sigma$ is separating in $M$. By Lemma~\ref{lem:closed}, $M$ is closed. Let $\bar U_1$ and $\bar U_2$ be the closures of the connected components of the complement of $\Sigma$ in $M$. Suppose first that $\dim M = 2$. Then $\Sigma$ is a circle, and since the time-$1$ map of a reversing diffeotopy interchanges $\bar U_1$ and $\bar U_2$, the classification of compact surfaces with boundary shows that $\bar U_1$ and $\bar U_2$ are either two copies of a closed orientable surface of genus $g\geq 0$ with a disk removed or two copies of a connected sum of $k \geq 1$ real projective planes with a disk removed. Since the first integral homology group of such an orientable surface has rank $2g$ and since the boundary circle maps to twice a class in homology in the non-orientable case, Lemma~\ref{lem:surjection_on_homology} implies that $\bar U_1$ and $\bar U_2$ are disks. Consequently, $M$ is a $2$-sphere. 

Suppose next that $\dim M \geq 3$. In this case $M$ is simply connected by Lemma~\ref{lem:surjection_on_pi1}, in particular orientable.
Lemma~\ref{lem:surjection_on_homology} moreover implies that $H_i(\bar U_1;\Z)=H_i(\bar U_2;\Z)=0$ for $i=1,\ldots,n-1$. The Mayer--Vietoris sequence then gives $H_i(M;\Z)=0$ for $i=1,\ldots,n-1$. Hence, $M$ is homeomorphic to an $(n+1)$-sphere by Poincar\'e duality, the Hurewicz theorem, Whitehead's theorem, and the generalized Poincar\'e conjecture.
\end{proof}

\subsection{Proof of Theorem~\ref{thm:top_class}, (i): the two-sided case} \label{sub:orientable_case}

Consider a Zoll family of closed embedded two-sided hypersurfaces on an $(n+1)$-manifold $M$ with $n \geq 1$. In this case the fiber bundle $\nu: \PT M \rightarrow P$ is coorientable. By Subsection~\ref{sub:coor_and_lifts}, the hypersurfaces in the family are hence diffeomorphic to $\Sph^n$. As in Section~\ref{sec:Berger_conjec} we consider the lifted Legendrian fiber bundle $\nu^+:\ST M \rightarrow P^+ $ with connected fibers $\Sigma^+ \cong \Sigma$ and the corresponding two-fold covering $P^+  \rightarrow P$ of manifolds.

To show that $M$ is homeomorphic to an $(n+1)$-sphere, by Lemma~\ref{lem:sphere_fibers_implies_sphere_total_space} it suffices to show that the hypersurfaces of the family in $M$ are pointed isotopically reversible. To see this, we fix a hypersurface $\Sigma$ of the family and a point $p\in\Sigma$. Let $z\in\ST_pM$ be a conormal covector of $\Sigma$ at $p$ and choose a smooth path $c\colon[0,1]\to\ST_pM$ from $z$ to $\tau(z)$. Then $\gamma:=\nu^+\circ c$ is a smooth path in $P^+$, i.e.\ a one-parameter family of cooriented hypersurfaces in the Zoll family through $p$, starting and ending at $\Sigma$ with opposite coorientations.

We pull back the bundle $\nu^+\colon\ST M\to P^+$ along $\gamma$ to the smooth fiber bundle
\[
  \gamma^*\nu^+=\{(t,y)\in[0,1]\times\ST M \mid \gamma(t)=\nu^+(y)\}
\]
over $[0,1]$ with fiber $\Sigma$, and denote by $\hat\gamma\colon\gamma^*\nu^+\to\ST M$ the projection to the second factor. Since $[0,1]$ is contractible, this bundle is trivial, i.e.\ there exists a diffeomorphism $F\colon[0,1]\times\Sigma\to\gamma^*\nu^+$ with $\pi_{[0,1]}(F(t,s))=t$. Composing with $\hat\gamma$ and the projection $\mu^+\colon\ST M\to M$ yields a smooth isotopy
\[
  \psi:=\mu^+\circ\hat\gamma\circ F\colon[0,1]\times\Sigma\to M
\]
with $\psi(\{t\}\times\Sigma)=\Sigma_{\gamma(t)}$, where $\Sigma_{\gamma(t)}$ denotes the cooriented hypersurface corresponding to $\gamma(t)$, and we may choose $F$ such that $\psi(0,\cdot)$ is the inclusion $\Sigma\hookrightarrow M$. Since $\nu^+(c(t))=\gamma(t)$, the map $t\mapsto(t,c(t))$ is a smooth section of $\gamma^*\nu^+$. Hence the curve $s\colon[0,1]\to\Sigma$, $s(t):=\mathrm{pr}_\Sigma\big(F^{-1}(t,c(t))\big)$, is smooth. It starts at $p$ and satisfies $\psi(t,s(t))=p$ for all $t\in[0,1]$. Applying the isotopy extension theorem to this curve in $\Sigma$ and precomposing $\psi$ with the resulting diffeotopy of $\Sigma$, we can assume that $\psi(t,p)=p$ for all $t\in[0,1]$. By the isotopy extension theorem, $\psi$ extends to a diffeotopy $\Psi$ of $M$, which then fixes $p$.
Since the coorientation of $\Sigma_{\gamma(t)}$ obtained by pushing forward that of $\Sigma$ via $\Psi_t$ and the one given by $\gamma(t)$ agree for $t=0$, they agree for all $t$ by continuity. Since $\Sigma_{\gamma(1)}$ is $\Sigma$ with the opposite coorientation, $\Psi_1$ preserves $\Sigma$ and reverses the orientation of its normal bundle. Hence $\Sigma$ is pointed isotopically reversible.

We now turn to the parameter manifold $P$. For $n+1=2$, $P^+$ is homeomorphic to $\Sph^2$ by Proposition~\ref{prp:even_rat_hom_sphere}. For $n+1\geq 3$ the long exact homotopy sequence of the bundle $\mu^+$ gives $\pi_k(\ST M) = 0$ for
$1 \leq k \leq n-1$. Substituting into the corresponding sequence of $\nu^+$ gives $\pi_k(P^+)=0$ in the same range. Now Poincaré duality, the Hurewicz theorem, Whitehead's theorem and the generalized Poincar\'e conjecture imply that $P^+$ is homeomorphic to an $(n+1)$-sphere; cf.~\cite[Section~3]{Ma26}. To show that $P$ is a real projective space, we argue as in \cite[Proposition~3.5]{Ma26}. Since the family is simple and $\nu$ is coorientable, $\mu^+$ embeds the fibers of $\nu^+$ into $M$; see Section~\ref{sec:Berger_conjec}. Hence, applying Subsection~\ref{sub:double_Legendrian} to the double Legendrian fibration $(\mu^+,\nu^+)$ of $\ST M$, $\nu^+$ embeds the fiber $\ST_pM$ of $\mu^+$ into $P^+$, and by the generalized Schönflies theorem its image separates $P^+$ into two balls. The deck transformation $\tau$ leaves this sphere invariant and interchanges the two balls, since the image of the sphere in $P$ is the hypersurface $\Sigma_p$, which is one-sided because $\mu$ is not coorientable; see Subsection~\ref{sub:double_Legendrian}, with the roles of $\mu$ and $\nu$ interchanged. Alternatively, this follows from Brouwer’s fixed point theorem, since $\tau$ acts freely on $P^+$. Hence $P$ is obtained from the closure of one of these balls by identifying boundary points via $\tau$. Since $\nu^+$ maps $\ST_pM$, on which $\tau$ acts as the antipodal map, $\tau$-equivariantly onto the boundary, the action on the boundary is conjugate to the antipodal map, and hence P is a real projective space by the Alexander trick.

\subsection{Proof of Theorem~\ref{thm:top_class}, (ii): the one-sided case}\label{sub:non_or_case}

Consider now a Zoll family of closed one-sided embedded  hypersurfaces on an $(n+1)$-manifold $M$ with $n \geq 1$.  By Subsection~\ref{sub:coor_and_lifts}, the hypersurfaces in the family are diffeomorphic to $\RP^n$. For any such hypersurface $\Sigma$, we glue two copies of the metric completion $M_0$ of $M\setminus \Sigma$ with respect to the induced length metric on $M\setminus \Sigma$ of a complete auxiliary Riemannian metric on $M$ along their boundaries via the nontrivial deck transformation of the two-fold covering map $\partial M_0 \rightarrow \Sigma$. The  map from the resulting manifold $\tilde M$ to $M$ defined as the identity on the two copies of $M\setminus \Sigma$ and as the natural projection onto $\Sigma$ on the boundaries of these copies in $\tilde M$ is a two-fold covering map, so that $\tilde M$ inherits a smooth structure from $M$, and the preimage $\tilde \Sigma$ of $\Sigma$ in $\tilde M$ is an $n$-dimensional two-sided sphere. 

For any other hypersurface $\Sigma'$ of the family, an isotopy from $\Sigma$ to $\Sigma'$ through hypersurfaces of the family, composed with the double cover $\tilde\Sigma\to\Sigma$, lifts to an isotopy in $\tilde M$ from $\tilde\Sigma$ to the preimage of $\Sigma'$, and so the latter is also a two-sided sphere. Hence the Zoll family on $M$ lifts to a Zoll family of closed two-sided embedded hyperspheres on $\tilde M$ with the same parameter space $P$, and the claim about $P$ follows from (i); cf.\ \cite[Section~2.2]{Ma26}. By (i), $\tilde M$ is homeomorphic to an $(n+1)$-sphere. By construction, the preimage $\tilde\Sigma\cong\Sph^n$ of $\Sigma$ separates $\tilde M$ into two components, which are interchanged by the deck transformation, and the latter restricts to the deck transformation of the two-fold covering $\tilde\Sigma\to\Sigma\cong\RP^n$. The same argument as for the parameter space in (i) now shows that $M$ is homeomorphic to $\RP^{n+1}$; cf.\ \cite[Theorem~2.11]{Ma26}.

\section{Deformations for Besse families of minimal hypersurfaces}\label{sec:deformations}

The proof of the two-dimensional version of Theorem~\ref{thm:zoll_finite_isometry}, \cite[Theorem~1.4]{AG26}, is based on a general equivariant construction, \cite[Theorems~1.1 and~1.2]{AG26}, valid in all dimensions, and a construction of certain invariant odd functions on the $2$-sphere, \cite[Theorem~5.1]{AG26}, which these theorems take as input. The following lemma generalizes the statement of \cite[Theorem~5.1]{AG26} from dimension two to all dimensions. In contrast to \cite{AG26}, where the function is constructed with the help of the classification of finite subgroups of $\mathrm{O}(3)$, our argument is a short averaging construction.

Here, by a linear function on $\Sph^{n+1}$ we mean the restriction to $\Sph^{n+1}$ of a linear function on $\R^{n+2}$. Note that a linear function on $\Sph^{n+1}$ is trivial if it vanishes on a nonempty open set in $\Sph^{n+1}$. Hence, two linear functions on $\Sph^{n+1}$ coincide if they coincide on a nonempty open set in $\Sph^{n+1}$. Also note that the space of all linear functions on $\Sph^{n+1}$ is invariant under the natural action of the orthogonal group $\mathrm{O}(n+2)$ on $L^2(\Sph^{n+1})$.

\begin{lem}\label{lem:lemma} For $n\geq 1$ and any finite subgroup $G < \mathrm{O}(n+2)$ with $-\mathrm{id} \notin G$ there exists a smooth odd function $f:\Sph^{n+1} \to \R$ which is $L^2$-orthogonal to the linear functions on $\Sph^{n+1}$ and has stabilizer $G$ in $\mathrm{O}(n+2)$.
\end{lem}
\begin{proof} For a finite subgroup $G < \mathrm{O}(n+2)$ let $\hat G$ be the finite subgroup of $\mathrm{O}(n+2)$ generated by $G$ and $-\mathrm{id}$. It contains $G$ as a subgroup of index $2$. Let $\sigma: \hat G \rightarrow \{\pm 1\}$ be the unique nontrivial homomorphism which is trivial on $G$, and let $\Lambda \subset \Sph^{n+1}$ be a fundamental domain for the action of $\hat G$ on $\Sph^{n+1}$ (e.g.\ $\Lambda = \{y\in \Sph^{n+1} \mid d(x,y) \leq d(x,gy) \text{ for all } g \in \hat G \setminus \{\mathrm{id}\} \}$ for some $x\in \Sph^{n+1}$ with trivial $\hat G$ isotropy). Let $\rho$ be a nonnegative smooth function on $\Sph^{n+1}$ with tiny support in the interior of $\Lambda$ such that $\rho$ attains its maximum at a unique point $x_0$ and has a trivial stabilizer in $\mathrm{O}(n+2)$ (which is for instance the case if the third-order Taylor expansion of $\rho$ at its maximum point in geodesic coordinates is sufficiently asymmetric, since an isometry fixing this point is determined by its differential there). The function $h:\Sph^{n+1} \rightarrow \R$ defined by
\[
		h(x) = \sum_{g\in \hat G} \sigma(g) \rho(gx)
\]
is then smooth, zero on an open set of almost full measure in $\Sph^{n+1}$, odd and has stabilizer $G$ in $\mathrm{O}(n+2)$. Indeed, if $k\in\mathrm{O}(n+2)$ leaves $h$ invariant, then it preserves the set where $h$ attains its maximum, which is the orbit $Gx_0$. After composing with an element of $G$, we can assume that $k$ fixes $x_0$. Since $h$ coincides with $\rho$ on a ball around $x_0$ containing the support of $\rho$, which $k$ preserves, $k$ leaves $\rho$ invariant, so $k=\mathrm{id}$ by our choice of $\rho$.

Let $h=h_{\mathrm{lin}} + f$ be the $L^2$-orthogonal decomposition of $h$ into its linear part $h_{\mathrm{lin}}$ and the complement of $h_{\mathrm{lin}}$ which we call $f$. Both $h_{\mathrm{lin}}$ and $f$ are smooth, odd and $G$-invariant. In particular, if $G$ has no fixed points on $\Sph^{n+1}$, then $h_{\mathrm{lin}}$ is trivial and the claim holds for $f$. In general, we additionally observe that $f$ is nontrivial, as $h$ is not linear, and equal to the linear function $-h_{\mathrm{lin}}$ on an open set of almost full measure in $\Sph^{n+1}$. By the discussion above about linear functions this implies that if some $k\in \mathrm{O}(n+2)$ leaves $f$ invariant, then it also leaves  $h_{\mathrm{lin}}$ and then also $h$ invariant. Hence, in any case the function $f$ satisfies all desired properties.
\end{proof}

Given this lemma, Theorem~\ref{thm:zoll_finite_isometry} follows from \cite[Theorems~1.1 and~1.2]{AG26}. Here the resulting sphere $M$ is disjoint from the origin in $\R^{n+2}$, as the spheres provided by \cite[Theorems~1.1 and~1.2]{AG26} are radial graphs over the round unit sphere. Moreover, $M$ is invariant under the linear action of $G$, which acts on it by isometries preserving the Zoll family, and every isometry of $M$ is the restriction of an element of $G$; see \cite[Theorem~1.2 and Corollary~4.19]{AG26}.

Finally, we turn to the applications of Theorem~\ref{thm:zoll_finite_isometry}.

\begin{proof}[Proof of Theorem~\ref{thm:non_embed_odd}]
For even $n\geq 2$ there exist finite subgroups $G$ of $\mathrm{SO}(n+2)$ of arbitrarily high order that act freely on $\Sph^{n+1}$ and do not contain the antipodal map. For instance, every cyclic group of odd order acts diagonally by rotations on the spherical join $\Sph^1 * \cdots*\Sph^1 = \Sph^{2k-1}$ of $k\geq 2$ copies of the unit circle in this way. In fact, every such $G$ has odd order, as a free linear involution of $\Sph^{n+1}$ is $-\mathrm{id}$. Moreover, in this case also the action of $G$ on $\R^{n+2}\setminus \{0\}$ is free. For every such nontrivial group $G$, Theorem~\ref{thm:zoll_finite_isometry} provides an embedded $G$-invariant sphere $M$ in $\R^{n+2}\setminus\{0\}$ with a $G$-invariant Zoll family of embedded minimal hypersurfaces and isometry group $G$. Since the action of $G$ on $M$ is isometric and free, the quotient $M/G$ is a manifold and inherits a Riemannian metric from $M$. Moreover, since $G$ acts freely on $\PT M$ and $\nu$ is $G$-equivariant, $\nu$ descends to an orbifold fiber bundle $\PT(M/G)=\PT M/G\to P/G$ with connected fibers $\nu^{-1}(\sigma)/G_\sigma$, so that the Zoll family on $M$ descends to a Besse family on $M/G$. Since $G \subset \mathrm{SO}(n+2)$ preserves the bounded component of $\R^{n+2}\setminus M$, it also preserves the orientation of $M$, and so the quotient $M/G$ is orientable. Hence, by Proposition~\ref{Berger_odd_orientable}, this Besse family is in fact Zoll, i.e.\ all its hypersurfaces are regular, though not simple Zoll, i.e.\ not all of its hypersurfaces are embedded, e.g. by Theorem~\ref{thm:top_class}. Indeed, since $\pi_1(M/G) \cong G$ is nontrivial of odd order, $M/G$ is homeomorphic to neither $\Sph^{n+1}$ nor $\RP^{n+1}$. In fact, by Lemma~\ref{lem:simple_charact} none of its hypersurfaces is embedded.

It remains to show that $M/G$ has trivial isometry group. Since $M$ is simply connected, every isometry of $M/G$ lifts to an isometry of $M$. A lift of a nontrivial isometry is not contained in $G$, contradicting the fact that the isometry group of $M$ is $G$ by Theorem~\ref{thm:zoll_finite_isometry}.
\end{proof}

\begin{rem}\label{rem:non-Zoll}
A similar argument as above yields Corollary~\ref{cor:Besse_orbifold}. For $n\geq 1$, let $G<\mathrm{O}(n+2)$ be a finite subgroup with $-\mathrm{id}\notin G$ containing a nontrivial element $g$ whose fixed subspace in $\R^{n+2}$ has dimension at least two, e.g.\ the group of order two generated by a reflection. Theorem~\ref{thm:zoll_finite_isometry} provides an embedded $G$-invariant sphere $M$ in $\R^{n+2}\setminus\{0\}$ with a $G$-invariant Zoll family of embedded minimal hypersurfaces and isometry group $G$, which descends to a Besse family of minimal hypersurfaces on the orbifold $M/G$. Since $M$ is a radial graph over $\Sph^{n+1}$, the fixed point set of $g$ on $M$ has positive dimension, so $M/G$ has non-isolated singularities. Hence, as in the proof of Lemma~\ref{lem:orbifold_manifold}, $\ST(M/G)$ is not a manifold. Let $P$ denote the parameter space of the Besse family on $M/G$. Since there is a covering map $\ST(M/G) \to\PT P$, obtained by passing to the quotient by $G$ from the corresponding $G$-equivariant covering for the Zoll family on $M$ (see Subsection~\ref{sub:coor_and_lifts}), also $\PT P$ and hence $P$ cannot be a manifold. The quotient $M/G$ does not have constant curvature, as otherwise $M$, being compact and simply connected, would be isometric to a round sphere and hence have isometry group $\mathrm{O}(n+2)$, contradicting Theorem~\ref{thm:zoll_finite_isometry}.
\end{rem}

\end{document}